\documentclass[10pt]{article}%

\usepackage[show]{ed}
\usepackage{mathrsfs}
\DeclareMathAlphabet{\mathpzc}{OT1}{pzc}{m}{it}

\usepackage{color}
\usepackage{amsmath}
\usepackage{amsthm}
\usepackage{graphicx}
\usepackage{amsfonts}
\usepackage{amssymb}%
\usepackage{latexsym}
\usepackage{psfrag}
\usepackage{subfig}
\usepackage{tikz}
\usetikzlibrary{patterns}
\usepackage{accents}

\newtheorem{theorem}{Theorem}[section]
\newtheorem{corollary}[theorem]{Corollary}
\newtheorem{definition}[theorem]{Definition}

\newtheorem{lemma}[theorem]{Lemma}
\newtheorem*{lemma*}{Lemma}
\newtheorem{proposition}[theorem]{Proposition}

\newtheorem{assumption}[theorem]{Assumption}
\newtheorem{example}[theorem]{Example}
\newtheorem{remark}[theorem]{Remark}
\newtheorem*{remark*}{Remark}

\newcounter{mycount}

\numberwithin{equation}{section}
\numberwithin{figure}{section}
\numberwithin{table}{section}

\usepackage[a4paper]{geometry}
\newcommand{\beq}{\begin{equation}}
\newcommand{\eeq}{\end{equation}}
\newcommand{\beqs}{\begin{equation*}}
\newcommand{\eeqs}{\end{equation*}}
\newcommand{\bit}{\begin{itemize}}
\newcommand{\eit}{\end{itemize}}
\newcommand{\ben}{\begin{enumerate}}
\newcommand{\een}{\end{enumerate}}
\newcommand{\bal}{\begin{align}}
\newcommand{\eal}{\end{align}}
\newcommand{\bals}{\begin{align*}}
\newcommand{\eals}{\end{align*}}
\newcommand{\bse}{\begin{subequations}}
\newcommand{\ese}{\end{subequations}}
\newcommand{\bpr}{\begin{proposition}}
\newcommand{\epr}{\end{proposition}}
\newcommand{\bre}{\begin{remark}}
\newcommand{\ere}{\end{remark}}
\newcommand{\bpf}{\begin{proof}}
\newcommand{\epf}{\end{proof}}
\newcommand{\ble}{\begin{lemma}}
\newcommand{\ele}{\end{lemma}}
\newcommand{\bco}{\begin{corollary}}
\newcommand{\eco}{\end{corollary}}
\newcommand{\bex}{\begin{example}}
\newcommand{\eex}{\end{example}}
\newcommand{\bth}{\begin{theorem}}
\newcommand{\enth}{\end{theorem}}

\newcommand{\noi}{\noindent}

\newcommand{\cA}{{\cal A}}

\newcommand{\cF}{{\cal F}}
\newcommand{\cH}{{\cal H}}
\newcommand{\cI}{{\cal I}}
\newcommand{\cL}{{\cal L}}
\newcommand{\opL}{\cL}

\newcommand{\cS}{{\cal S}}
\newcommand{\cT}{{\cal T}}

\newcommand{\cM}{{\cal M}}

\newcommand{\bx}{x}%\mathbf{x}}

\newcommand{\ba}{a}%\mathbf{a}}

\newcommand{\bQ}{Q}%\mathbf{Q}}

\newcommand{\supp}{\mathrm{supp}}

\newcommand{\re}{{\rm e}}
\newcommand{\ri}{{\rm i}}
\newcommand{\rd}{{\rm d}}

\newcommand{\Rea}{\mathbb{R}}

\newcommand{\diam}{\mathop{{\rm diam}}}

\newcommand{\pdiff}[2]{\frac{\partial #1}{\partial #2}}

\newcommand{\OR}{{B_R}}

\newcommand{\GR}{{\partial B_R}}

\newcommand{\ngus}{|\nabla u|^2}
\newcommand{\nurs}{|u_r|^2}
\newcommand{\ngtus}{\ngus - \nurs}
\newcommand{\gu}{\nabla u}

\newcommand{\nvs}{|v|^2}
\newcommand{\ngvs}{|\nabla v|^2}

\newcommand{\nvrs}{|v_r|^2}

\newcommand{\gv}{\nabla v}

\newcommand{\gvb}{\overline{\nabla v}}

\newcommand{\vb}{\overline{v}}

\newcommand{\LtG}{{L^2(\partial B_R)}}

\newcommand{\HoDkk}{{H^1_{k}(B_R)}}

\newcommand{\tendi}{\rightarrow \infty}
\newcommand{\tendo}{\rightarrow 0}

\newcommand{\abs}[1]{|#1|}

\newcommand{\NN}{\mathbb{N}}

\newcommand{\Hilb}{\cH}

\def\XXint#1#2#3{{\setbox0=\hbox{$#1{#2#3}{\int}$}
     \vcenter{\hbox{$#2#3$}}\kern-.5\wd0}}

\newcommand*{\N}[1]{\left\|#1\right\|}
\newcommand*{\norm}[1]{\N{#1}}

\newcommand{\vertiii}[1]{{\left\vert\kern-0.25ex\left\vert\kern-0.25ex\left\vert #1 
    \right\vert\kern-0.25ex\right\vert\kern-0.25ex\right\vert}}

\usepackage{hyperref}
\definecolor{myblue}{rgb}{0,0,0.6}
\hypersetup{colorlinks=true,
linkcolor=myblue,citecolor=myblue,filecolor=myblue,urlcolor=myblue}

\allowdisplaybreaks[4]
\newcommand{\ton}{\text{ on }}
\newcommand{\tin}{\text{ in }}
\newcommand{\tfa}{\text{ for all }}

\newcommand{\tfor}{\text{ for }}
\newcommand{\tas}{\text{ as }}
\newcommand{\tand}{\text{ and }}
\newcommand{\tst}{\text{ such that }}

\newcommand{\hatx}{\widehat{\bx}}

\newcommand{\domain}{B}

\definecolor{escol}{rgb}{0,0,0.8}
\definecolor{estcol}{rgb}{0,0.8,0}

\newcommand{\DtN}{{\rm DtN}_k}
\newcommand{\CDtN}{{C_{\rm DtN}}}

\newcommand{\Cqo}{{C_{\rm qo}}}
\newcommand{\Ccont}{{C_{\rm cont}}}
\newcommand{\Csol}{{C_{\rm sol}}}

\newcommand{\CFEM}{{C_{\rm approx}}}
\newcommand{\CHt}{{C_{{\rm split}, H^2}}}

\newcommand{\CA}{{C_{{\rm split},\cA}}}

\newcommand{\hFEM}{h}

\newcommand{\uhigh}{u_{H^2}}
\newcommand{\ulow}{u_{\mathcal A}}

\newcommand{\Pilow}{\Pi_{L}}%\flat}}
\newcommand{\Pihigh}{\Pi_{H}}%\sharp}}
\newcommand{\hsc}{{k^{-1}}}

\newcommand{\pa}{\partial}

\newcommand{\Schwartz}{\mathscr{S}}
\newcommand{\Schwartzdual}{\Schwartz^*}

\newcommand{\mythmname}[1]{\textbf{\emph{(#1.)}}}

\newcommand{\trace}{}

\begin{document}

\title{A simple proof that the $hp$-FEM does not suffer from the pollution effect for the constant-coefficient full-space Helmholtz equation
}

\author{
E. A.~Spence\footnotemark[1]
}

\date{\today}

\renewcommand{\thefootnote}{\fnsymbol{footnote}}

\footnotetext[1]{Department of Mathematical Sciences, University of Bath, Bath, BA2 7AY, UK, \tt E.A.Spence@bath.ac.uk}

\renewcommand{\thefootnote}{\arabic{footnote}}

\maketitle

\begin{abstract}
In $d$ dimensions, approximating an arbitrary function oscillating with frequency $\lesssim k$ requires $\sim k^d$ degrees of freedom. A numerical method for solving the Helmholtz equation (with wavenumber $k$) suffers from the pollution effect if, as $k\tendi$, the total number of degrees of freedom needed to maintain accuracy grows faster than this natural threshold. 

While the $h$-version of the finite element method (FEM) (where accuracy is increased by decreasing the meshwidth $h$ and keeping the polynomial degree $p$ fixed) suffers from the pollution effect, the celebrated papers \cite{MeSa:10, MeSa:11,EsMe:12,MePaSa:13} showed that the $hp$-FEM 
(where accuracy is increased by decreasing the meshwidth $h$ and increasing the polynomial degree $p$)
applied to a variety of constant-coefficient Helmholtz problems does not suffer from the pollution effect.

The heart of the proofs of these results is a PDE result splitting the solution of the Helmholtz equation into ``high'' and ``low'' frequency components. 
In this expository paper we prove this splitting for the constant-coefficient  Helmholtz equation in 
full space (i.e., in $\Rea^d$) using \emph{only} integration by parts and elementary properties of the Fourier transform; this is in contrast to the proof for this set-up in \cite{MeSa:10} which uses somewhat-involved bounds on Bessel and Hankel functions. The proof in this paper is motivated by the recent proof in \cite{LSW3} of this splitting for the variable-coefficient Helmholtz equation in full space; indeed, the proof in \cite{LSW3} uses more-sophisticated tools that reduce to the elementary ones above for constant coefficients.

We combine this splitting with (i) standard arguments about convergence of the FEM applied to the Helmholtz equation (the so-called ``Schatz argument'', which we reproduce here) and (ii) polynomial-approximation results (which we quote from the literature without proof) to give a simple proof that the $hp$-FEM does not suffer from the pollution effect for the constant-coefficient full-space Helmholtz equation.
\end{abstract}

\section{Introduction and motivation}

When computing approximations with a numerical method to the solution of the Helmholtz equation
\beq\label{eq:Helmholtz_intro}
\Delta u +k^2 u=0, 
\eeq
with wavenumber $k>0$, a fundamental question is:
\begin{quotation}
\noi How must the number of degrees of freedom increase with $k$ to maintain accuracy as $k\tendi$?
\end{quotation}
In this paper, we consider the finite-element method (FEM) using piecewise-polynomial approximation spaces, with meshwidth $h$ and polynomial degree $p$. Recall that in the $h$-version of the FEM, accuracy is increased by decreasing $h$, in the $p$-version accuracy is increased by increasing $p$, and in the $hp$-version accuracy is increased by \emph{both} decreasing $h$ \emph{and} increasing $p$. In the context of these methods, the question above is 
\begin{quotation}
\noi How quickly must $h$ decrease with $k$ and/or $p$ increase with $k$ to maintain accuracy as $k\tendi$?
\end{quotation}
This question has been the subject of sustained interest since the late 1980's, with initially answers obtained for 1-d problems 
\cite{AzKeSt:88, IhBa:95, IhBa:97} (see also \cite[Chapter 4]{Ih:98}), and now answers obtained for 2- and 3-d problems in general geometries, both for ``standard'' FEMs  \cite{Me:95,  MeSa:10, MeSa:11, EsMe:12, ZhWu:13, Wu:14, DuWu:15, LiWu:19, LSW3, GaChNiTo:22, GS3} 
and for variations of these, e.g. discontinuous Galerkin methods
\cite{FeWu:09, FeWu:11, MePaSa:13, ZhWu:21} and multiscale methods \cite{GaPe:15, Ch:16, BaChGo:17, Pe:17, OhVe:18, CaWu:20, ChVa:20}.
Moreover, there is large current interest in this question when the Helmholtz equation \eqref{eq:Helmholtz_intro} is replaced by its variable-coefficient generalisation $\nabla\cdot(A\nabla u ) +k^2 nu=0$ \cite{BrGaPe:17, ChNi:20, GaSpWu:20, GrSa:20, LSW4, MaAlSc:21, LSW3, GLSW1, BeChMe:22} 
or even the time-harmonic Maxwell equations \cite{NiTo:20, MeSa:21, MeSa:22}.

A highlight of this body of research is the result from 
\cite{MeSa:10, MeSa:11, EsMe:12, MePaSa:13} that the $hp$-FEM does not suffer from the \emph{pollution effect}; i.e., accuracy can be maintained with a choice of the number of degrees of freedom growing like $k^d$, where ``accuracy'' here means that the computed solution is \emph{quasi-optimal}; see \eqref{eq:qo_gen} below. This is contrast to the $h$-version of the FEM which, e.g., with $p=1$, needs the total number of degrees of freedom to grow like $k^{2d}$ to maintain accuracy in this sense. Having the number of degrees of freedom growing like $k^d$ is the natural threshold for this problem since an oscillatory function in $d$ dimensions with frequency $\lesssim k$
requires $\sim k^d$ degrees of freedom to be well-approximated by piecewise polynomials;
this is expected in 1-d from the Nyquist--Shannon--Whittaker sampling theorem \cite{Wh:15, Sh:49} (see, e.g., \cite[Theorem 5.21.1]{BaNaBe:00}) and from the recent results in general dimension in \cite{G1}.

The proofs that the $hp$-FEM does not suffer from the pollution effect in \cite{MeSa:10, MeSa:11, EsMe:12, MePaSa:13}  
consist of the following three ingredients.
\ben
\item Sufficient conditions for FEM solutions to be quasi-optimal
originating from the ideas of Schatz \cite{Sc:74}  (related to the classic ``Aubin--Nitsche trick"), and then developed by Sauter \cite{Sa:06}.
This argument is now well-known; e.g., the essence of it appears in the books \cite[\S5.7]{BrSc:08}, \cite[\S4.4.2]{Ih:98}.
\item Results from \cite[Appendices B and C]{MeSa:10} about how well the $hp$-FEM spaces approximate analytic functions  (these are more-sophisticated versions of the standard piecewise-polynomial approximation results appearing in, e.g., \cite[Chapter 4]{BrSc:08}, \cite[Chapter 17]{Ci:91}).
\item A PDE result splitting the solution of the Helmholtz equation into ``high'' and ``low'' frequency components (this was the heart of the paper \cite{MeSa:10}).
\een
The motivation for the present paper was the realisation that, for the constant-coefficient Helmholtz equation \eqref{eq:Helmholtz_intro} posed in $\Rea^d$, 
given a bound on the solution in terms of the data (which can be proved using essentially only integration by parts),
the splitting in Point 3 above can be proved using \emph{only} elementary properties of the Fourier transform, thus making the key ideas behind this active area of research accessible to a wide audience.
The proof in this paper is motivated by the recent proof in \cite{LSW3} of this splitting for the \emph{variable-coefficient} Helmholtz equation in $\Rea^d$; the proof in \cite{LSW3} uses more-sophisticated tools that reduce to the elementary ones above in the constant-coefficient case -- we discuss this further in Section \ref{sec:LSW}.

\vspace{-1ex}

\paragraph{Plan of the paper.}

After recalling basic facts about the Helmholtz equation (and why it is important) in Section \ref{sec:Helmholtz}, 
Sections \ref{sec:Schatz}, \ref{sec:hpFEM}, and  \ref{sec:splitting} concern Points 1, 2, and 3 above, respectively.
Section \ref{sec:splitting} states the splitting for the constant-coefficient full-space Helmholtz equation, and then uses the results of Points 1-3 to prove that the $hp$-FEM does not suffer from the pollution effect when applied to this problem.
Section \ref{sec:Fourier} recaps the basic properties of the Fourier transform, and then Section \ref{sec:proof} proves the splitting stated in Section \ref{sec:splitting} using  the material in Section \ref{sec:Fourier}.

\vspace{-1ex}

\paragraph{Prerequisite knowledge required.}

This paper is pitched at a reader who knows about the following two topics.
\bit
\item The analysis of the FEM applied to Poisson's equation $\Delta u=-f$ at the level of, for example, \cite[Chapters 1-5]{BrSc:08} (covering weak derivatives, the Sobolev space $H^1$, Green's identity, the Lax--Milgram theorem, C\'ea's lemma, construction of finite-element spaces, and basic polynomial-approximation theory).
\item Basic properties of the Fourier transform (the inversion formula, the Fourier transform of derivatives, and Plancherel's theorem).  
\eit
The paper is largely self-contained; the main exceptions are the following.
\bit 
\item Some of the proofs of the background results about the Helmholtz equation in Section \ref{sec:Helmholtz} are only sketched, with references given to full proofs in the books \cite{CoKr:83, KiHe:15, Ne:01}.
\item The $hp$-approximation results in Section \ref{sec:hpFEM} are given without proof -- we hope  that the reader can take these on trust, following discussion comparing and contrasting these with more-standard polynomial-approximation results appearing in, e.g., \cite{Ci:91, BrSc:08}.
\eit

\section{The Helmholtz equation}\label{sec:Helmholtz}

\subsection{Where does the Helmholtz equation come from?}

The Helmholtz equation \eqref{eq:Helmholtz_intro}
is arguably the simplest possible model of wave propagation.
For example, if we look for solutions of the wave equation
\beq\label{eq:wave}
\pdiff{^2 U}{t^2} - c^2 \Delta U =0
\eeq
in the form 
\beq\label{eq:time_harmonic}
U(\bx,t) =
u(\bx) \re^{\pm \ri \omega t},
\eeq
then the function $u(\bx)$ satisfies the Helmholtz equation \eqref{eq:Helmholtz_intro} with $k = \omega/c$. Assuming a similar dependence on time reduces the Maxwell equations to the so-called time-harmonic Maxwell equations, which, in certain situations, can be further reduced to the Helmholtz equation (see, e.g., \cite[\S1.4.3]{Ih:98} \cite[Remark 2.1]{MoSp:19}). Similarly, the time-harmonic elastic wave equation (often called the Navier equation) also reduces to the Helmholtz equation 
in certain situations (see, e.g., \cite[\S1.2]{Ih:98}).

\subsection{The model Helmholtz problem considered in this paper}

Given $f\in L^2(\Rea^d)$ with compact support, let $u$ be the solution of 
\beq\label{eq:Helmholtz}
k^{-2}\Delta u  + u = -f \,\,\text{ in } \Rea^d, \,\,d=2,3,
\eeq
satisfying the \emph{Sommerfeld radiation condition}
\beq\label{eq:src}
k^{-1}\pdiff{u}{r}(\bx) - \ri  u(\bx) = o \left( \frac{1}{r^{(d-1)/2}}\right)
\eeq
as $r:= |\bx|\tendi$, uniformly in $\hatx:= \bx/r$. We make two immediate remarks.

(i) We have multiplied the Helmholtz equation \eqref{eq:Helmholtz_intro} by $k^{-2}$ (and added a source term); we see below 
how this rescaling by $k^{-2}$ allows us to keep better track of the $k$-dependence.

(ii) We show below (Theorem \ref{thm:Fred}) that the solution to \eqref{eq:Helmholtz}-\eqref{eq:src} exists and is unique. In fact, the solution can be written down explicitly as an integral of $f$ against the fundamental solution of the Helmholtz equation (see \eqref{eq:convolution} below); however the results in this paper (both about the Helmholtz equation itself and the FEM applied to it)
 hold for much more general Helmholtz problems where the solution cannot be written down explicitly.

\subsection{The meaning of the radiation condition \eqref{eq:src}}

The Sommerfeld radiation condition \eqref{eq:src} expresses mathematically that, with the choice $\re^{-\ri \omega t}$ in \eqref{eq:time_harmonic}, the scattered wave moves away from the obstacle towards infinity. Indeed, one can show (see Theorem \ref{thm:AtWi} below) that \eqref{eq:src} implies that 
\beq\label{eq:AtWi_intro}
u(\bx) = \frac{\re^{\ri k r}}{r^{(d-1)/2}}\left( F_0(x/r ) + O\left(\frac{1}{r}\right)\right) \tas r\to \infty,
\eeq
for some smooth function $F_0(x/r)$ (i.e., a function of the non-radial variables). 
Recalling that $k=\omega/c$, we find that the corresponding solution $U(\bx,t) := 
 u(\bx) \re^{-\ri \omega t}$ 
 of the wave equation \eqref{eq:wave}
satisfies
\beqs
U(\bx,t) 
= \frac{\re^{\ri k (r-ct)}}{r^{(d-1)/2}}\left( F_0( x/r ) + O\left(\frac{1}{r}\right)\right),
\eeqs
and we recognise $\re^{\ri k (r-ct)}$ as a wave travelling in the positive $r$ direction (i.e., away from the support of $f$ -- the source of the waves) as time $t$ increases.
\footnote{If we seek solutions of the wave equation as $U(x,t)= u(x)\re^{\ri \omega t}$, then the radiation condition corresponding to outgoing waves is 
$\partial_r u +\ri ku = o \big(r^{(1-d)/2}\big)$ as $r\to \infty$.}
From now on, we therefore say that a solution of the Helmholtz equation satisfying \eqref{eq:src} is \emph{outgoing}.

\subsection{The variational formulation of the Helmholtz equation}\label{sec:vf}

\paragraph{The goal.}
Find a Hilbert space $\cH$, a sesquilinear form $a(\cdot,\cdot)$ (i.e., linear in the first argument and antilinear in the second argument), and an antilinear functional $F(\cdot)$ such that 
the Helmholtz problem \eqref{eq:Helmholtz}-\eqref{eq:src} is equivalent to the variational problem:
\beq\label{eq:EDPvar2}
\text{ find } u \in \cH\tst\,\, a(u,v)=F(v) \,\, \tfa v\in  \cH.
\eeq
\paragraph{Motivation.} The FEM is a special case of the Galerkin method, which is a method for obtaining approximations to the solutions of variational problems of the form \eqref{eq:EDPvar2}; see \S\ref{sec:Galerkin} below.

\vspace{-1ex}

\paragraph{What is the issue?}
The natural space in which to seek weak solutions of second-order linear elliptic PDEs (such as the Helmholtz equation \eqref{eq:Helmholtz}) posed on a domain $D$ is 
\beqs
H^1(D):= \big\{ v\in L^2(D) \text{ such that } \nabla v \in L^2(D)\big\}.
\eeqs
The Helmholtz problem \eqref{eq:Helmholtz}-\eqref{eq:src} is posed on the unbounded domain $\Rea^d$, but the radiation condition \eqref{eq:src} (and its consequence the expansion 
\eqref{eq:AtWi_intro}) implies that the solution $u$ is in neither $L^2(\Rea^d)$ nor $H^1(\Rea^d)$. It is, however, in $H^1(B_R)$ for every $R>0$, where $B_R$ is the ball of radius $R$ centred at the origin; we say that such a function is in $H^1_{\rm loc}(\Rea^d)$.

\vspace{-1ex}

\paragraph{The solution.}
The summary of how to deal with this issue is as follows.
\bit
\item[Step 1:] Prove that the solution to \eqref{eq:Helmholtz}-\eqref{eq:src}, if it exists, is unique. 
\item[Step 2:] Consider the Helmholtz equation $k^{-2}\Delta u+u=0$ in the exterior of $B_R$ and satisfying the radiation condition \eqref{eq:src}. 
Use separation of variables to find an explicit expression for the solution; the uniqueness proof from Step 1 also shows that this solution is unique. 
\item[Step 3:] Use the explicit solution in the exterior of $B_R$ from Step 2 and Green's identity to create a variational problem on $B_R$ (with $\cH= H^1(B_R)$) whose solution is $u|_{B_R}$; the uniqueness result in Step 1 implies that the solution of this variational problem is unique.
\item[Step 4:] Prove wellposedness of the variational problem on $B_R$ (which proves existence of the solution of \eqref{eq:Helmholtz}-\eqref{eq:src}).
\eit

\paragraph{Which results in this subsection are needed for the analysis of the Galerkin method in \S\ref{sec:Schatz}?}
The analysis of the Galerkin method in \S\ref{sec:Schatz} uses (i) the definition of the variational problem on $B_R$ \eqref{eq:EDPa}, (ii) the fact that the solution of this variational problem exists and is unique (i.e., the result of Theorem \ref{thm:Fred}), and (iii) the two properties of the sesquilinear form in Lemma \ref{lem:propa}. We highlight this explicitly for a reader who wants to get to the Galerkin-method analysis as quickly as possible and is willing to take on trust these PDE results.

We now go through each of Steps 1-4 above. 

\vspace{-1ex}

\paragraph{Step 1 (Uniqueness of the solution).}

\begin{theorem}\mythmname{Atkinson--Wilcox expansion}\label{thm:AtWi}
If $u\in H^1_{\rm loc}(\Rea^d\setminus \overline{B_{R_0}})$ is an outgoing solution of  $k^{-2}\Delta u+u=0$ in $\Rea^d\setminus \overline{B_{R_0}}$ for some $R_0>0$, then there exist smooth functions $F_n$ such that, for any $R_1>R_0$, 
\beq\label{eq:AtWi}
u(x) =\frac{\re^{\ri k r}}{r^{(d-1)/2}}\sum_{n=0}^\infty\frac{F_n(x/r)}{r^n} \quad\tfor r:=|x|\geq R_1
\eeq
(so that, in particular, the expansion \eqref{eq:AtWi_intro} holds),
where the sum in \eqref{eq:AtWi} (and all its derivatives) converges absolutely and uniformly.
Furthermore, if $F_0\equiv 0$ then $u\equiv0$ in $\Rea^d\setminus \overline{B_{R_1}}$.
\end{theorem}

\bpf[Sketch proof]
This proof uses Green's integral representation for $u$ outside $B_{R_1}$, i.e.,
\beq\label{eq:GreenIR}
u(x) = \int_{\partial B_{R_1}}\left(u(y) \pdiff{\Phi_k(x,y)}{\nu(y)} - \pdiff{u}{\nu}(y) \Phi_k(x,y)\right) \rd s(y)\quad \tfor |x|>R_1,
\eeq
where $\nu(x)=x/r$ (i.e., $\nu$ is the outward-pointing unit normal vector to $\partial B_{R_1}$)
and  $\Phi_k(x,y)$ is the fundamental solution (satisfying $(\Delta_y +k^2)\Phi_k(x,y)= -\delta(x-y)$) defined by  
\beq\label{eq:fund}
\Phi_k(x,y):= 
\displaystyle{\frac{\ri }{4}H_0^{(1)}\big(k|x-y|\big)},  \quad d=2,\quad
:=\displaystyle{\frac{\re^{\ri k |x-y|}}{4\pi |x-y|}}, \quad d=3.
\eeq
The proof of \eqref{eq:GreenIR} can be found in, e.g., \cite[Theorem 3.3]{CoKr:83}; this proof 
holds when $u\in C^2(\Rea^d\setminus \overline{B_{R_1}})$, but this is true by elliptic regularity. Indeed, 
if $u\in H^1_{\rm loc}(\Rea^d\setminus \overline{B_{R_0}})$ is a weak solution of the Helmholtz equation, then applying the elliptic regularity bound of Corollary \ref{cor:elliptic_regularity} 
to a cut-off function multiplied by $u$ 
(as in the proof of Corollary \ref{cor:H2} below)
implies that $u\in H^2(\Rea^d\setminus \overline{B_{R_1}})$ for any $R_1>R_0$. Repeating this process with $u$ replaced by its second derivatives, we find that 
$u\in H^4(\Rea^d\setminus \overline{B_{R_1}})$, and, by induction, $u \in H^n(\Rea^d\setminus \overline{B_{R_1}})$ for all $n$. The Sobolev embedding theorem (see, e.g., \cite[Theorem 3.26]{Mc:00}) then implies that $u\in C^\infty(\Rea^d\setminus \overline{B_{R_1}})$. 

Using the large $|x|$ asymptotics of $\Phi_k(x,y)$ (obtained by direct calculation), we obtain the expansion \eqref{eq:AtWi} from \eqref{eq:GreenIR}; see \cite[Theorem 3.6]{CoKr:83}. The requirement that $u$ satisfies the Helmholtz equation implies a recurrence relation for the $F_n$ (see \cite[Corollary 3.8]{CoKr:83}), which implies that if $F_0\equiv 0$, then $F_n\equiv 0$ for all $n$, and thus $u(x)=0$ for $|x|\geq R_1$. % (for every $R_1>R_0$). 
%Finally, to show that $u(x)=0$ everywhere, observe that Green's integral representation implies that $u$ is real analytic in $\Rea^d$ (i.e., can be expanded as a convergent power series about every point in $\Rea^d$); thus if $u(x)=0$ for $|x|\geq R_1$, then $u\equiv 0$.
%Get expansion of $u$ from expansion of fundamental solution. Proof shows in fact that $u= \re^{\ri k r} r^{-(d-1)/2} \sum_{n=0}^\infty F_n(x/r) r^{-n}$. If $F_0$ then all the $F_n$ are zero and $u$ is zero.
\epf

\begin{theorem}\mythmname{Rellich's uniqueness theorem}\label{thm:unique}
The solution of \eqref{eq:Helmholtz}-\eqref{eq:src} % in $H^1_{\rm loc}(\Rea^d)$
 is unique.
\end{theorem}

\bpf[Sketch proof]
We need to show that if $f\equiv 0$, then $u\equiv 0$.
Applying Green's identity 
\beq\label{eq:Green}
\langle \partial_\nu u, v \rangle_{\partial D} = \int_D \nabla u\cdot \overline{\nabla v } + \overline{v} \Delta u
\eeq
(where $\langle\cdot, \cdot\rangle_{\partial D}$ in this situation is just $(\cdot,\cdot)_{L^2(\partial D)}$) with $D=B_R$ and $v=u$, and then taking the imaginary part, we find that
$\Im \int_{\partial B_R} \overline{u} \,\partial u/\partial r=0$
(this application of Green's identity is justified since $u\in C^\infty(\Rea^d)$ using the reasoning in the proof of Theorem \ref{thm:AtWi}).
The expansion \eqref{eq:AtWi_intro} implies that, as $R\to \infty$, $\Im \int_{\partial B_R} \overline{u} \,\partial u/\partial r$ tends to the integral over the unit sphere of $|F_0(x/r)|^2$; therefore $F_0\equiv 0$. By Theorem \ref{thm:AtWi}, $u(x)=0$ for $|x|\geq R_1$
 for every $R_1>0$, and thus $u\equiv 0$ in $\Rea^d$.
\epf

\paragraph{Step 2 (Explicit expression for solution of outgoing Helmholtz in exterior of ball).}

For simplicity we just consider $d=2$. The expressions for $d=3$ are similar, with spherical harmonics replacing the trigonometric polynomials $\re^{\ri n \theta}$; see, e.g., 
\cite[\S2.6]{Ne:01}, 
\cite[\S3]{ChMo:08}, \cite[\S3]{MeSa:10}.
%\cite[Equations 3.5 and 3.6]{ChMo:08} \cite[\S2.6.3]{Ne:01}, \cite[Equations 3.7 and 3.10]{MeSa:10}. 

\begin{theorem}\mythmname{Explicit solution in exterior of ball}
Let $d=2$. Given $g\in H^{1/2}(\partial B_R)$, let 
\beq\label{eq:sov}
v(r,\theta) := \frac{1}{2\pi} \sum_{n=-\infty}^\infty \frac{H^{(1)}_n(kr)}{ H^{(1)}_n(kR)} \re^{\ri n \theta} \widehat{g}(n),
\quad\text{ where }\quad
\widehat{g}(n):= \int_0^{2\pi} \re^{-\ri n\theta} g(R, \theta)\, \rd \theta.
\eeq
Then $v$ is the unique outgoing solution to 
\beq\label{eq:DtN1}
(-k^{-2}\Delta - 1)v=0 \quad\tin \Rea^d \setminus \overline{B_R} \quad\tand\quad \trace v =g \ton \partial B_R.
\eeq
\end{theorem}

\bpf[Sketch proof]
The proof of uniqueness in Theorem \ref{thm:unique} for the solution of \eqref{eq:Helmholtz}-\eqref{eq:src} also shows that the outgoing solution of \eqref{eq:DtN1} is unique (the only change is that we now apply Green's identity
in $B_{R_1}\setminus B_R$ and send $R_1\to \infty$; the integral over $\partial B_R$ is zero thanks to the Dirichlet boundary condition).
By the definition of the Hankel function $H_n^{(1)}(z)$, the product $H_n^{(1)}(kr) \re^{\ri n \theta}$ satisfies the Helmholtz equation. The boundary condition on $\partial B_R$ is satisfied by the inversion theorem for Fourier series.
The asymptotics of $H_n^{(1)}(z)$ for fixed $n$ as $z\to \infty$ (see, e.g., \cite[Equation 10.17.5]{Di:22}) show that, for fixed $n$, $H_n^{(1)}(kr)$ is outgoing; these asymptotics, however, do not hold uniformly with respect to $n$, which makes proving that 
the sum in \eqref{eq:sov} is outgoing more difficult. One option is to use the uniform asymptotics of Hankel functions (see, e.g., \cite[\S10.20]{Di:22}); an alternative proof (using Green's integral representation theorem) is given in \cite[Theorem 2.37]{KiHe:15} (this proof is for $d=3$, but the ideas are exactly the same for $d=2$).
\epf

We now define the \emph{Dirichlet-to-Neumann map in the exterior of $B_R$}, 
$\DtN: H^{1/2}(\GR)\to H^{-1/2}(\GR)$, by, given $g\in H^{1/2}(\partial B_R)$, with $v$ is given by \eqref{eq:sov},
\beq\label{eq:DtN}
\DtN g := k^{-1} \partial_r v,
\quad\text{ i.e., }\quad\DtN g (\theta)= \frac{1}{2\pi} \sum_{n=-\infty}^\infty \frac{H^{(1)'}_n(kR)}{ H^{(1)}_n(kR)} \exp(\ri n \theta) \widehat{g}(n).
\eeq
For the analogous expression when $d=3$, see, 
e.g., \cite[Equations 3.5 and 3.6]{ChMo:08} \cite[\S2.6.3]{Ne:01}, \cite[Equations 3.7 and 3.10]{MeSa:10}. 

\paragraph{Step 3 (The Helmholtz variational formulation on $B_R$.)} In this step we use the fact that 
the traces (i.e., restrictions to the boundary) of functions in $H^1(B_R)$ live in the space $H^{1/2}(\partial B_R)$, and their normal derivatives (when they exist) live in the space $H^{-1/2}(\partial B_R)$; however, 
further information about these fractional Sobolev spaces is not needed in the rest of the paper.

Let $R>0$ be large enough so that $\supp f \subset B_R$. The variational formulation of \eqref{eq:Helmholtz}-\eqref{eq:src} is then 
\beq\label{eq:EDPvar}
\text{ find } \widetilde{u} \in H^1(B_R) \tst \,\, a(\widetilde{u},v)=F(v) \,\, \tfa v\in H^1(B_R),
\eeq
where
\beq\label{eq:EDPa}
a(\widetilde{u},v):= \int_{B_R} 
\Big(k^{-2}\nabla \widetilde{u}\cdot\gvb
 - \widetilde{u}\vb\Big) - k^{-1}\big\langle \DtN \trace \widetilde{u},\trace v\big\rangle_{\partial B_R}\quad\tand\quad
F(v):= \int_{B_R} f\, \vb.
\eeq
The pairing $\langle \cdot,\cdot\rangle_{\GR}$ in $a(\cdot,\cdot)$ is the duality pairing between $H^{-1/2}(\partial B_R)$ and $H^{1/2}(\partial B_R)$; however, for the reader unfamiliar with this concept, there is little harm in just thinking of it as the $\LtG$ inner product $(\cdot,\cdot)_{\LtG}$.

\ble\mythmname{Equivalence of the formulations}\label{lem:equiv}
If $u$ is a solution of \eqref{eq:Helmholtz}-\eqref{eq:src}, then $u|_{B_R}$ is a solution of the variational problem \eqref{eq:EDPvar}. 
Conversely, if $\widetilde{u}$ is a solution of this variational problem, then there exists a solution $u$ of \eqref{eq:Helmholtz}-\eqref{eq:src} such that $u|_{B_R}= \widetilde{u}$.
\ele

\bpf[Sketch proof]
Recall that Green's identity states that, if $v\in H^1(D)$ and $u \in H^1(D)$ with $\Delta u \in L^2(D)$, then \eqref{eq:Green} holds.
Given a solution of  \eqref{eq:Helmholtz}-\eqref{eq:src}, apply \eqref{eq:Green} in $B_R$; the definitions of $\DtN$ \eqref{eq:DtN} and $a(\cdot,\cdot)$ \eqref{eq:EDPa} then imply that $\widetilde{u}$ satisfies \eqref{eq:EDPvar}. Conversely, given $\widetilde{u}$ a solution of \eqref{eq:EDPvar}, let $u:= \widetilde{u}$ in $B_R$ and $u:= v$ in $(B_R)^c$ where $v$ is defined by \eqref{eq:DtN1}/\eqref{eq:sov} with $g:= \widetilde{u}|_{\partial B_R}$. 
Green's identity \eqref{eq:Green} and the definitions of $\DtN$ and $a(\cdot,\cdot)$ then imply that $u$ is a (weak) solution of \eqref{eq:Helmholtz} in $\Rea^d$.
\epf

Combining Theorem \ref{thm:unique} and Lemma \ref{lem:equiv} we obtain the following result.

\begin{corollary}\mythmname{Uniqueness of the solution of the variational problem}\label{cor:uniqueness}
If $a(u,v)= 0$ for all $v\in H^1(B_R)$ then $u=0$.
\end{corollary}

\paragraph{Step 4 (Wellposedness of the variational problem).}

We work with the weighted norm $\|\cdot\|_{H^1_k(B_R)}$ defined by 
\beq\label{eq:1knorm}
\|u\|^2_{H^1_k(B_R)}:= 
\|k^{-1}\nabla u \|^2_{L^2(B_R)} + \|u \|^2_{L^2(B_R)};
\eeq
this is a special case of the weighted norm
\beq\label{eq:weighted_norms}
\N{v}^2_{H^m_k(B_R)}:= 
\sum_{0\leq |\alpha|\leq m}
\N{(k^{-1}\partial)^\alpha v}^2_{L^2(B_R)}.
\eeq
The rationale for using these norms is that if a function $v$ oscillates with frequency $k$, then we expect $|(k^{-1}\partial)^\alpha v|\sim  |v|$ for all $\alpha$; this is true, e.g., if $v(\bx) = \exp(\ri k \bx\cdot\ba)$. 
\footnote{
%\interfootnotelinepenalty=100000
Many papers on numerical analysis of the Helmholtz equation use the weighted $H^1$ norm 
\beqs
\|u\|^2_{H^1_k(B_R)}:= \|\nabla u \|^2_{L^2(B_R)} + k^2\|u \|^2_{L^2(B_R)};
\eeqs
we use \eqref{eq:1knorm}/\eqref{eq:weighted_norms} instead since weighting the $j$th derivative by $k^{-j}$ is easier to keep track of than weighting it by $k^{-j+1}$ (especially for high derivatives).}

\ble\mythmname{Properties of $a(\cdot,\cdot)$}\label{lem:propa}

(i) (Continuity) 
Given $k_0, R_0>0$ there exists $\Ccont>0$ such that for all $k\geq k_0$ and $R\geq R_0$, 
\beq\label{eq:continuity}
|a(u,v)|\leq \Ccont \N{u}_{H^1_k(B_R)} \N{v}_{H^1_k(B_R)} \quad\tfa u, v \in H^1(B_R).
\eeq

(ii) (G\aa rding inequality) 
\begin{align}\label{eq:Garding}
\Re a(v,v) &\geq  \N{v}^2_{H^1_k(\domain_R)} - 2\N{v}^2_{L^2(\domain_R)} \quad \tfa v\in H^1(B_R).
\end{align}
\ele

Lemma \ref{lem:propa} is proved using the following two properties of $\DtN$.

\ble\mythmname{Key properties of $\DtN$}\label{lem:DtN}

(i) Given $k_0,R_0>0$ there exists $\CDtN_1= \CDtN_1(k_0 R_0)$ such that for all $k\geq k_0$ and $R\geq R_0$,
\beq\label{eq:CDtN1}
\big|\big\langle \DtN\trace u, \trace v\rangle_{\partial B_R}\big\rangle\big| \leq k\,\CDtN_1 \N{u}_{H^1_k(\OR)}  \N{v}_{H^1_k(\OR)} \quad\tfa u,v \in H^1(\OR).
\eeq

(ii) 
\beq\label{eq:CDtN2}
- \Re \big\langle \DtN \phi,\phi\big\rangle_{\GR} \geq 0
 \quad\tfa \phi\in H^{1/2}(\GR).
\eeq
\ele

\bpf[References for the proof of Lemma \ref{lem:DtN}]
Both (i) and (ii) are proved using the expression for $\DtN$ in terms of Bessel and Hankel functions (i.e., \eqref{eq:DtN} when $d=2$) in \cite[Lemma 3.3]{MeSa:10} (see also \cite[Theorem 2.6.4]{Ne:01} and \cite[Lemma 2.1]{ChMo:08} for \eqref{eq:CDtN2}).
\epf

\bpf[Proof of Lemma \ref{lem:propa}]
(i) follows from the definition of $a(\cdot,\cdot)$ \eqref{eq:EDPa}, the Cauchy-Schwarz inequality, the definition of $\|\cdot\|_{H^1_k(\OR)}$ \eqref{eq:1knorm}, and the inequality \eqref{eq:CDtN1}.

(ii) follows from the definitions of $a(\cdot,\cdot)$ and $\|\cdot\|_{H^1_k(\OR)}$ 
and the inequality \eqref{eq:CDtN2}.
\epf

\bre\mythmname{$a(\cdot,\cdot)$ is not coercive on $H^1(B_R)$ for $k$ sufficiently large}
The G\aa rding inequality \eqref{eq:Garding} shows that, up to the lower-order $L^2$ term,  $a(\cdot,\cdot)$ is coercive on $H^1_k(B_R)$. 
We highlight that,  for sufficiently large $k$, $a(\cdot,\cdot)$ is \emph{not} coercive, i.e., there does not exist $\gamma>0$ such that $|a(v,v)|\geq \gamma\|v\|_{H^1_k(B_R)}^2$ for all $v \in H^1(B_R)$, and thus neither the Lax--Milgram lemma nor C\'ea's lemma are applicable. Indeed, if $\lambda_j$ is a Dirichlet eigenvalue of $-\Delta$ in $B_R$ with eigenfunction $u_j$, then $a(u_j,u_j)=0$ when $k^2=\lambda_j$. A little more work shows that if $k^2\geq \lambda_1$ (i.e., the smallest Dirichlet eigenvalue in $B_R$) then there exists $v\in H^1(B_R)$ such that $a(v,v)=0$, and thus $a(\cdot,\cdot)$ is not  coercive for $k^2\geq \lambda_1$; see, e.g., \cite[Lemma 6.34]{Sp:15}.
\ere

\bth\mythmname{Wellposedness of the Helmholtz problem}\label{thm:Fred}
The variational problem \eqref{eq:EDPvar} has a unique solution which depends continuously on the data $f$.
\end{theorem}

\bpf
By Lemma \ref{lem:equiv}, it is sufficient to prove that the variational problem \eqref{eq:EDPvar} has a unique solution which depends continuously on the data.
Since the sesquilinear form is continuous \eqref{eq:continuity} and satisfies the G\aa rding inequality \eqref{eq:Garding}, Fredholm theory implies 
that existence of a solution to the variational problem and continuous dependence of the solution on the data both follow from uniqueness; see, e.g., \cite[\S6.2.8]{Ev:98}, \cite[Theorem 2.34]{Mc:00}.
\epf

\bre\mythmname{Approximating $\DtN$}
Implementing the operator $\DtN$ appearing in $a(\cdot,\cdot)$ \eqref{eq:EDPa} is computationally expensive, and so in practice one seeks to approximate this operator by \emph{either} imposing an absorbing boundary condition on $\GR$, \emph{or} using a perfectly-matched layer (PML), \emph{or} using boundary integral equations (so-called ``FEM-BEM coupling''); see, e.g., the overview in \cite[\S3]{Ih:98}.
For simplicity, in this paper we analyse the FEM assuming that $\DtN$ is realised exactly; recent $k$-explicit results on the error incurred (on the PDE level) by approximating $\DtN$ by absorbing boundary conditions or PML can be found in \cite{GaLaSp:21} and \cite{GaLaSp:21a}, respectively.
\ere

\subsection{The $k$-dependence of the Helmholtz solution operator}

The wellposedness result of Theorem \ref{thm:Fred} implies that 
\beq\label{eq:Csol}
\N{u}_{H^1_k(\OR)}\leq \Csol \N{f}_{L^2(\OR)}
\eeq
for some $\Csol<\infty$. A natural question is then:~how does $\Csol$ depend on $k$ and $R$?
We define more precisely
\beqs
\Csol(k,R):= 
\sup_{f\in L^2(B_R),\, \|f\|_{L^2(B_R)} = 1}
\|u\|_{H^1_k(B_R)};
\eeqs
i.e., $\Csol(k,R)$ is the operator norm of the map $L^2(\OR)\ni f \mapsto u\in H^1(\OR)$, where $u$ is the outgoing solution of the Helmholtz equation \eqref{eq:Helmholtz}. 

The large-$k$ behaviour of solutions of the Helmholtz equation is dictated by geometric-optic rays; therefore, in general, proving $k$-explicit bounds on the Helmholtz solution operator requires analysis that takes into account the behaviour of these rays (see, e.g., \cite{DyZw:19}). However, for the simple model problem \eqref{eq:Helmholtz}-\eqref{eq:src}, the following bound on $\Csol$ can be obtained by multiplying the Helmholtz equation \eqref{eq:Helmholtz} by a judiciously-chosen test function and integrating by parts. 

\begin{theorem}\mythmname{Morawetz bound on $\Csol$}\label{thm:Morawetz}
For all $k>0$ and $R>0$,
\beq\label{eq:Morawetz_bound}
\Csol \leq 2 kR \sqrt{ 1 + \left(\frac{d-1}{2kR}\right)^2}.
\eeq
\end{theorem}
This bound was essentially proved in \cite{MoLu:68, Mo:75} (although the bound does not quite appear in this form in those papers); see also \cite[Lemma 3.5]{ChMo:08} and \cite[Equations 1.9 and 1.10]{GrPeSp:19}.
The details of the proof of this bound are not needed in the rest of the paper, but for completeness (and since they involve essentially only integrating by parts) we include them in Appendix \ref{app:Morawetz}.

The bound \eqref{eq:Morawetz_bound} is sharp in its $kR$ dependence for $kR$ large; indeed, by considering $u(x)= \re^{\ri k x_1} \chi(|x|/R)$ for $\chi \in C^\infty$ supported in $[0,1)$ one can show that given $k_0,R_0>0$ there exists $C>0$ such that $\Csol \geq C kR$ for all $k\geq k_0$ and $R\geq R_0$.

The convergence theory of the FEM relies on knowing how smooth the solution of the underlying PDE is (see \S\ref{sec:hpFEM} below).
When the data is in $L^2$, the solution of a second-order linear elliptic PDE is in $H^2$, provided the coefficients and boundary have sufficient regularity. For the problem \eqref{eq:Helmholtz}-\eqref{eq:src}, the following corollary shows that this $H^2$ regularity is achieved.

\begin{corollary}\mythmname{Bound on $H^2_k$ norm of the Helmholtz solution}\label{cor:H2}
Given $k_0, R_0>0$, there exists $C>0$ such that the solution $u$ of \eqref{eq:Helmholtz}-\eqref{eq:src} with $\supp f \subset B_R$ with $R\geq R_0$ satisfies
\beq\label{eq:H2}
\N{u}_{H^2_k(B_R)} \leq C kR \N{f}_{L^2(B_R)} \quad \tfa k\geq k_0. 
\eeq
\end{corollary}

\bpf[Sketch proof]
Corollary \ref{cor:elliptic_regularity} below uses the Fourier transform to prove $H^2$ regularity of the solution of $(-k^{-2}\Delta+1)v=g \in L^2(\Rea^d)$. To apply this bound to 
the solution of \eqref{eq:Helmholtz}-\eqref{eq:src}, let $\varphi \in C^\infty_{\rm comp}(\Rea^d, [0,1])$ be equal to one on $B_{R}$ and vanish outside $B_{2R}$, and then let $v:= \varphi u$ so that $g=(-k^{-2}\Delta+1)(\varphi u)=(-k^{-2}\Delta-1)(\varphi u)+ 2 \varphi u$. The bound \eqref{eq:H2} then follows by bounding $\|g\|_{L^2(B_{2R})}$ in terms of $\|f\|_{L^2(B_R)}$ using the PDE \eqref{eq:Helmholtz} and the bounds \eqref{eq:Csol} and \eqref{eq:Morawetz_bound}.
\epf

\section{The Galerkin method, quasioptimality, and the Schatz argument
}\label{sec:Schatz}

\subsection{Definition of the Galerkin method and quasioptimality}\label{sec:Galerkin}

In the setting of \S\ref{sec:vf}, 
let $\cH_N$ be a finite-dimensional subspace of $\cH$ (which is $H^1(B_R)$ for our Helmholtz problem).
The Galerkin method applied to the variational problem \eqref{eq:EDPvar2} is 
\beq\label{eq:FEM}
\text{ find } u_N \in \cH_N\tst\,\, a(u_N,v_N)=F(v_N) \,\, \tfa v_N\in  \cH_N.
\eeq
The FEM is the Galerkin method \eqref{eq:FEM} with $\cH_N$ consisting of piecewise polynomials (we describe in \S\ref{sec:hpFEM} the specific assumptions we make on $\cH_N$).

Given a sequence of finite-dimensional spaces $\{\cH_N\}$, a standard error bound one seeks to prove on the sequence of Galerkin solutions $\{u_N\}$ is the following:~there exists a $\Cqo>0$ and $N_0\in \NN$ such that, for $N\geq N_0$,
\beq\label{eq:qo_gen}
\N{u-u_N}_{\cH} \leq \Cqo \min_{v_N\in \cH_N}\N{u-v_N}_{\cH}.
\eeq
If such a bound holds then we say that the sequence of Galerkin solutions (or, more informally, the Galerkin method itself) is 
\emph{quasioptimal}. (If \eqref{eq:qo_gen} holds with $\Cqo=1$ then the method is \emph{optimal}.)

\subsection{The adjoint solution operator with $L^2$ data}

The Schatz argument crucially relies on properties of the adjoint solution operator.

\begin{definition}\mythmname{Adjoint solution operator $\cS^*$}
Given $f\in L^2(\OR)$, let $\cS^*f  \in H^1(\OR)$ be defined by 
\beq\label{eq:S*vp}
a(v, \cS^*f) = (v,f)_{L^2(\OR)} \,\,\tfa v\in H^1(\OR).
\eeq
\end{definition}

The following lemma shows that our knowledge about outgoing Helmholtz solutions immediately gives us knowledge about $\cS^*$.

\ble\label{lem:Helmholtz_adjoint}
If $\cS^*$ is defined as in \eqref{eq:S*vp} then
\beqs%\label{eq:S*fkey}
a(\overline{\cS^*f}, v)= (\overline{f},v)_{L^2(\OR)}\quad\tfa v\in H^1(\OR).
\eeqs
i.e.,~$\cS^*f$ is the complex-conjugate of the outgoing Helmholtz solution with data $\overline{f}$.
\ele

\bpf[Sketch proof]
Green's identity \eqref{eq:Green} and the radiation condition \eqref{eq:src} show that $\big\langle \DtN \psi, \overline{\phi}\big\rangle_{\GR} =\big\langle \DtN \phi, \overline{\psi}\big\rangle_{\GR} $ for all $\phi,\psi\in H^{1/2}(\GR)$.
This implies that $a(\overline{v},u)=a(\overline{u},v)$ for all $u,v$, which implies the result.
\epf

\subsection{The Schatz argument:~sufficient conditions for quasioptimality}

\ble\mythmname{Sufficient conditions for quasi-optimality}\label{lem:Schatz}
If 
\beq\label{eq:etasuff}
\eta(\cH_N): = \sup_{0 \neq f\in L^2(\OR)} \min_{v_N\in\cH_N} \frac{\N{\cS^*f- v_N}_{\HoDkk}}{\N{f}_{L^2(\OR)}}
\leq \frac{1}{2\Ccont},
\eeq
then the Galerkin solution $u_N$ to the variational problem \eqref{eq:FEM} exists, is unique, and satisfies the bound
\beq\label{eq:qo}
\N{u-u_N}_{H^1_k(\OR)} \leq 2\Ccont \left(\min_{v_N\in V_N} \N{u-v_N}_{H^1_k(\OR)}\right).
\eeq
\ele

Observe that the quantity $\eta(\cH_N)$ in \eqref{eq:etasuff} measures how well solutions of the adjoint problem are approximated in the space $\cH_N$.

\bpf
We first prove the result under the assumption that the Galerkin solution $u_N$ exists.

\paragraph{Step 1:~Use the G\aa rding inequality and Galerkin orthogonality.}

Setting $v=v_N$ in \eqref{eq:EDPvar2} and subtracting this from \eqref{eq:FEM}, we obtain the \emph{Galerkin orthogonality} that
\beq\label{eq:GO}
a(u-u_N,v_N) =0 \quad\tfa v_N\in\cH_N.
\eeq
Using (in this order) the G\aa rding inequality \eqref{eq:Garding}, Galerkin orthogonality \eqref{eq:GO}, and continuity of $a(\cdot,\cdot)$ \eqref{eq:continuity}, we have that, for any $v_N\in \cH_N$,
\begin{align}\nonumber
\N{u-u_N}^2_{\HoDkk} 
&\leq \Re a(u-u_N, u-u_N) +  2\N{u-u_N}^2_{L^2(\OR)}\\ \nonumber
&= \Re a(u-u_N, u-v_N) +  2\N{u-u_N}^2_{L^2(\OR)}\\
& \leq \Ccont \N{u-u_N}_{\HoDkk} \N{u-v_N}_{\HoDkk}
+  2\N{u-u_N}^2_{L^2(\OR)}.\label{eq:key1}
\end{align}

\paragraph{Step 2:~Prove an Aubin-Nitsche type bound using the definition of $\eta(\cH_N)$.}

The quasioptimal error bound \eqref{eq:qo} under the condition \eqref{eq:etasuff} follows from \eqref{eq:key1} if we can prove that 
\beq\label{eq:Schatz1}
\N{u-u_N}_{L^2(\OR)} \leq \Ccont \eta(\cH_N) \N{u-u_N}_{\HoDkk}.
\eeq
By the definition of $S^*$ \eqref{eq:S*vp}, Galerkin orthogonality \eqref{eq:GO}, and continuity \eqref{eq:continuity},
\begin{align}\nonumber
\N{u-u_N}^2_{L^2(\OR)} = a\big(u-u_N, S^*(u-u_N)\big) &= a\big(u-u_N, S^*(u-u_N)-v_N \big) \\
&\hspace{-1.5cm} \leq \Ccont \N{u-u_N}_{\HoDkk} \N{S^*(u-u_N)-v_N}_{\HoDkk}\label{eq:S2}
\end{align}
for any $v_N\in \Hilb_N$.
The definition of $\eta(\Hilb_N)$ \eqref{eq:etasuff} implies that there exists a $w_{h,p}\in \Hilb_N$ such that
\beqs
\N{S^*(u-u_N)-w_{h,p}}_{\HoDkk} \leq \eta(\Hilb_N) \N{u-u_N}_{L^2(\OR)}.
\eeqs
Using this last inequality in \eqref{eq:S2}, we obtain \eqref{eq:Schatz1}.

\paragraph{Step 3:~Prove that $u_N$ exists.}

We have so far assumed that $u_N$ exists. Recall that an $N\times N$ matrix $B$ is invertible if and only if  $B$ has full rank, which is the case if and only if the only solution of $Bx=0$ is $x=0$. Therefore, to show that $u_N$ exists, we only need to show that $u_N$ is unique. Seeking a contradiction, suppose that there exists a $\widetilde{u}_N\in \Hilb_N$ such that 
$a(\widetilde{u}_N,v_N)=0 \tfa v_N\in \Hilb_N.$
Let $\widetilde{u}$ be such that 
\beq\label{eq:S4}
a(\widetilde{u},v)=0 \quad\tfa v\in H^1(B_R);
\eeq
thus $\widetilde{u}_N$ is the Galerkin approximation to $\widetilde{u}$.
Repeating the argument in the first part of the proof we see that if \eqref{eq:etasuff} holds then the quasi-optimal error bound \eqref{eq:qo} holds (with $u$ replaced by $\widetilde{u}$ and $u_N$ replaced by $\widetilde{u}_N$). By Corollary \ref{cor:uniqueness}, the only solution to the variational problem \eqref{eq:S4} is $\widetilde{u}=0$, and then \eqref{eq:qo} implies that $\widetilde{u}_N=0$.
We have therefore shown that the solution $u_N$ exists under the condition \eqref{eq:etasuff} and the proof is complete.
\epf

\bre\mythmname{Bibliographic remarks}
We describe \eqref{eq:Schatz1}  as an ``Aubin-Nitsche-type bound'', since the argument that obtains \eqref{eq:Schatz1} 
was first introduced in the %%%%self-adjoint, -- Nitsche does self-adjoint, Aubin not
coercive case by Aubin \cite[Theorem 3.1]{Au:67} and Nitsche \cite{Ni:68} (see, e.g., \cite[Theorem 19.1]{Ci:91} for \eqref{eq:Schatz1} stated as the ``Aubin--Nitsche lemma").
Schatz \cite{Sc:74} considered second-order linear elliptic PDEs satisfying a G\aa rding inequality (such as \eqref{eq:Garding})
proving existence and uniqueness of the Galerkin solution for $N$ sufficiently large. The fact that these arguments also give quasioptimality was recognised in \cite[Theorem 3.1]{AzKeSt:88}.
The concept of $\eta(\cH_N)$ and emphasis on the role of ``adjoint approximability'' are due to Sauter \cite{Sa:06}.
\ere

\section{Recap of approximation results in $hp$-FEM spaces}\label{sec:hpFEM}

The result that the $hp$-FEM does not suffer from the pollution effect is proved under the following two assumptions on the finite-dimensional subspaces; % (Assumptions \ref{ass:FEM1} and \ref{ass:FEM2} below).
these assumptions describe how well (as a function of $h$ and $p$) the spaces approximate functions with a given regularity.

We highlight immediately that both these assumptions are satisfied by $hp$-finite-element spaces with curved elements that fit $\partial B_R$ exactly, provided that the triangulations are quasi-uniform 
and are constructed by refining a fixed triangulation that has analytic element maps (see Theorem \ref{thm:assFEM1} below and the associated discussion).
Nevertheless, we formulate these properties as specific assumptions to make it clear the actual properties of the subspaces that are needed in the proof of the result that the $hp$-FEM does not suffer from the pollution effect. Since we are ultimately thinking of the subspaces in these assumptions as $hp$-finite-element spaces, we denote the sequence of these subspaces as $\{ \cH_{h,p}\}_{h>0,p\in\mathbb{Z}^+}$.

\begin{assumption}
\mythmname{Approximation of functions with finite regularity} 
\label{ass:FEM1}
Let $\{ \cH_{h,p}\}_{h>0,p\in\mathbb{Z}^+}$ be a sequence of finite-dimensional subspaces of $H^1(B_R)$.
Given $s,d$ with $d$ the spatial dimension and $s>d/2$, there exists $\CFEM_1>0$ 
such that if $v \in H^s(B_R)$ and $p\geq s-1$, then
\beq\label{eq:assFEM1}
\min_{w_{h,p} \in \cH_{h,p}} \N{v-w_{h,p}}_{H^1_k(B_R)} \leq \CFEM_1 \left(\frac{hk}{p}\right)^{s-1} \left(1 + \frac{hk}{p}\right) \N{v}_{H^s_k(B_R)}.
\eeq
%where $|v|_{H^s(B_R)}$ is the $H^s(B_R)$ semi-norm; i.e. $|v|_{H^s(B_R)}^2:= \sum_{|\alpha|=s}\|\partial^{\alpha} v\|_{L^2(B_R)}^2$.
\end{assumption}

\paragraph{Discussion of Assumption \ref{ass:FEM1}.}
%Linking the bound \eqref{eq:assFEM1} to more-standard polynomial approximation results.}
A standard polynomial approximation result is the following. For $d=2,3$, given $s\geq 2$ and $p\geq s-1$, there exists $C>0$ 
such that if $v \in H^s(D)$ and $m=0$ or $1$, then
\beq\label{eq:polyapprox}
\N{v-\cI^h v }_{H^m(D)} \leq C h^{s-m} |v|_{H^s(D)},
\eeq
where $\cI^h$ is a global interpolation operator; see, e.g., \cite[Equation 4.4.28]{BrSc:08}, \cite[Theorem 17.1]{Ci:91}. The approximation result \eqref{eq:assFEM1} is a generalisation of \eqref{eq:polyapprox} which (i) makes explicit the dependence on $p$ of the constant $C$ in \eqref{eq:polyapprox}, and (ii) works in norms weighted with $k$.

\paragraph{Motivation for Assumption \ref{ass:FEM2}.}
Assumption \ref{ass:FEM1} is about approximating  in $k$-weighted norms an arbitrary function in $H^s$ for some $s>0$; note that the constant $\CFEM_1$ depends on $s$ in an unspecified way, and so we cannot use the bound \eqref{eq:assFEM1} for arbitrarily-large $s$, and hence also arbitrarily-large $p$, since $p$ is tied to $s$ via $p\geq s-1$. 
The next assumption, Assumption \ref{ass:FEM2}, allows us to take arbitrarily-large $p$; the price one pays is that the function being approximated must be analytic. 

Before stating Assumption \ref{ass:FEM2}, we recall the relationship between derivative bounds and analyticity for families of functions depending on $k$.

\ble[$k$-explicit analyticity]\label{lem:analytic}
Let $D$ be a bounded open subset of $\Rea^d$ and 
let $u\in C^\infty(D)$ 
be a family of functions depending on $k$.

(i) If 
there exist $C, C_u>0$, independent of $\alpha$, such that 
\beq\label{eq:analytic1}
\N{\partial^\alpha u}_{L^2(D)}\leq C_u (C k)^{\vert \alpha\vert }\quad\text{ for all multiindices } \alpha,
\eeq
then $u$ is real analytic in $D$ and its power series has infinite radius of convergence, i.e., $u$ can be extended to an entire function on $\Rea^d$.

(ii) If 
there exist $C, C_u>0$, independent of $\alpha$, such that 
\beqs
\N{\partial^\alpha u}_{L^2(D)}\leq C_u (Ck)^{\vert \alpha\vert  } \vert \alpha\vert  ! \quad\text{ for all multiindices } \alpha,
\eeqs
then $u$ is real analytic in $D$ with radius of convergence of its power series proportional to $(C k)^{-1}$. 

(iii) If 
there exist $C, C_u>0$, independent of $\alpha$, such that 
\beqs
\N{\partial^\alpha u}_{L^2(D)}\leq C_u C^{\vert \alpha\vert  }\max \big\{ \vert \alpha\vert  , k \big\}^{\vert \alpha\vert  } \quad\text{ for all multiindices } \alpha,
\eeqs
then $u$ is real analytic in $D$ with radius of convergence of its power series proportional to $C^{-1}$ and independent of $k$.
\ele

\bpf[Sketch proof]
In each case, use the Sobolev embedding theorem (see, e.g., \cite[Theorem 3.26]{Mc:00}) to obtain a bound on $\|\partial^\alpha u\|_{L^\infty(D)}$, and then use this to bound the Lagrange form of the remainder in the Taylor series; see, e.g., 
\cite[Proof of Lemma C.2]{MeSa:10}. 
\epf

In the rest of the paper, we only use the class of functions in Part (i) of Lemma \ref{lem:analytic}, but the classes in Parts (ii) and (iii) are included for context.

\begin{assumption}[Approximation 
in the finite-dimensional subspace of the class of functions in Part (i) of Lemma \ref{lem:analytic}]
\label{ass:FEM2}
Let  $\{\cH_{h,p}\}_{h>0,p\in\mathbb{Z}^+}$ be a sequence of finite-dimensional subspaces of $H^1(B_R)$, and 
suppose $v\in C^\infty(B_R)$ is such that, given $k_0>0$ there exists $C_1,C_2>0$ such that
\beq\label{eq:better_than_analytic}
\N{(k^{-1}\partial)^\alpha v}_{L^2(B_R)}\leq C_1 (C_2)^{\vert \alpha\vert  } \quad\tfa \alpha \text{ and for all } k\geq k_0.
\eeq
Given $\widetilde{C}$, there exist $\sigma,\CFEM_2>0$, depending on $C_2$ and $\widetilde{C}$ (but not $C_1$), such that, if $k\geq k_0$ and $k, h$, and $p$ satisfy
\beq\label{eq:assFEM2cond}
\frac{h}{R} + \frac{hk}{p}\leq \widetilde{C},
\eeq
then
\beq\label{eq:assFEM2}
\min_{w_{h,p} \in \cH_{h,p}} \N{v-w_{h,p}}_{H^1_k(B_R)} \leq C_1 \CFEM_2 
\left[
\frac{1}{kR}\left(\frac{h/R}{h/R+\sigma}\right)^p
 + \left(\frac{hk}{\sigma p}\right)^p 
\right].
\eeq
\end{assumption}

\paragraph{Discussion of Assumption \ref{ass:FEM2}.}
The key points about the bound \eqref{eq:assFEM2} are the following.
\bit
\item[(i)] For fixed $p$ and $k$, as $h\to 0$, the right hand side of \eqref{eq:assFEM2} is $O( h^p)$ (just like the right-hand side of \eqref{eq:assFEM1} with $s=p+1$).
\item[(ii)] If $hk/(\sigma p)<1$ then the right-hand side of \eqref{eq:assFEM2} decreases exponentially as $p\tendi$.
\item[(iii)] The quantities $h/R$, $hk/p$, and $kR$ are dimensionless, and thus the right-hand sides of \eqref{eq:assFEM2cond} and \eqref{eq:assFEM2} involve only dimensionless quantities.
\eit

\paragraph{Approximation spaces satisfying Assumptions \ref{ass:FEM1} and \ref{ass:FEM2}.}

Let 
$(\cT_{\hFEM})_{0<h\leq h_0}$ (with $h$ the maximum element diameter) be a sequence of triangulations of $B_R$, with each element $K\in\cT_{\hFEM}$ the image of a reference element $\widehat{K}$ (a reference triangle in 2-d and a reference tetrahedron in 3-d) under the map $F_K : \widehat{K}\to K$. As is standard, we assume there are no hanging nodes and that the element maps of elements sharing an edge or face induce the same parametrisation on that edge or face. 
We consider the $hp$-finite-element spaces
\beq\label{eq:cHh}
 \cH^{p}(\cT_h):= \Big\{v \in H^1(\OR)
 : \text{for each } K \in \cT_{\hFEM}, \,v|_K \circ F_K \text{ is a polynomial of degree $\leq p$} 
 \Big\}.
 \eeq
Since $B_R$ is curved, we consider triangulations with curved elements that fit $\partial B_R$ exactly (thus avoiding the issue of analysing the non-conforming error coming from using simplicial triangulations; see, e.g., \cite[Chapter 10]{BrSc:08}).  
Recall that the family $(\cT_h)_{0<h\leq h_0}$ is \emph{quasi-uniform} if there exists $C>0$ such that 
\beqs
h:=\max_{K \in \cT_h} \diam(K) \leq C \min_{K \in \cT_h} \diam(K) \quad\tfa 0<h\leq h_0;
\eeqs
for such triangulations, the dimension of $ \cH^{p}(\cT_h)$ is proportional to $(p/(h/R))^d$.

%\begin{theorem}\mythmname{Conditions under which Assumptions \ref{ass:FEM1} and \ref{ass:FEM2} hold}\label{thm:assFEM1}
%Assume that $d=2,3,$ and $(\cT_h)_{0<h\leq h_0}$ is a family of quasi-uniform triangulations. Assume further that each $K \in \cT$ is the image of a triangle/tetrahedron under the image of a bi-Lipschitz map (i.e., both the map and its inverse are Lipschitz). Then $\cH^{p}(\cT_h)$ defined by \eqref{eq:cHh} satisfies Assumption \ref{ass:FEM1}. 
%\end{theorem}

%\bpf[References for the proof]
%This follows from \cite[Theorem B.4]{MeSa:10} (a result about approximation on the reference element) and a scaling argument (see \cite[Bottom of Page 1895]{MeSa:10}). The result  \cite[Theorem B.4]{MeSa:10} builds on the results of \cite[Theorem 4.1]{Mu:97}, \cite{Zh:08}, and \cite{GuZh:09}; see the discussion at the start of \cite[Appendix B]{MeSa:10}.
%\epf

\begin{theorem}\mythmname{Conditions under which Assumptions \ref{ass:FEM1} and \ref{ass:FEM2} hold}\label{thm:assFEM1}
If $(\cT_h)_{0<h\leq h_0}$ satisfies \cite[Assumption 5.2]{MeSa:10}, then  $\cH^{p}(\cT_h)$ defined by \eqref{eq:cHh} satisfies Assumptions \ref{ass:FEM1} and \ref{ass:FEM2}.
\end{theorem}

Informally, \cite[Assumption 5.2]{MeSa:10} is that $(\cT_h)_{0<h\leq h_0}$ is quasi-uniform with each element map $F_K$ the composition of a affine map and an analytic map; \cite[Remark 5.2]{MeSa:11} notes that $(\cT_h)_{0<h\leq h_0}$ 
satisfying this assumption can be constructed by refining a fixed triangulation that has analytic element maps.

\bpf[References for the proof of Theorem \ref{thm:assFEM1}]
That Assumption \ref{ass:FEM1} holds follows from \cite[Theorem B.4]{MeSa:10} (a result about approximation on the reference element) and a scaling argument (see \cite[Bottom of Page 1895]{MeSa:10}). We note that the result  \cite[Theorem B.4]{MeSa:10} builds on the results of \cite[Theorem 4.1]{Mu:97}, \cite{Zh:08}, and \cite{GuZh:09}; see the discussion at the start of \cite[Appendix B]{MeSa:10}.

For Assumption \ref{ass:FEM2}, the bound \eqref{eq:assFEM2} is proved in the course of  \cite[Proof of Theorem 5.5]{MeSa:10}, see the last equation on \cite[Page 1896]{MeSa:10}; note that 
(i) we have simplified this equation using the assumption \eqref{eq:assFEM2cond}, and (ii) the weighted $H^1$ norm in \cite{MeSa:10} is $k$ times $\|\cdot\|_{H^1_k(B_R)}$ defined by \eqref{eq:1knorm}.
\epf

\section{The splitting of the Helmholtz solution and the proof that the $hp$-FEM does not suffer from the pollution effect}\label{sec:splitting}

\subsection{Statement of the splitting}

The crucial result used to prove that the $hp$-FEM applied to the problem \eqref{eq:Helmholtz}-\eqref{eq:src} does not suffer from the pollution effect is the following (with the proof contained in \S\ref{sec:proof}).

\begin{theorem}\mythmname{Splitting of the Helmholtz solution}\label{thm:splitting}
Given $k_0, R_0>0$, there exists $\CHt$, $\CA>0$ such that the following holds.
Given $f\in L^2(B_R)$ with $R\geq R_0$, let $u$ satisfy the Helmholtz equation \eqref{eq:Helmholtz}
 and the Sommerfeld radiation condition \eqref{eq:src}.
Then
 \beq\label{eq:splitting}
u|_{B_R}= \uhigh+ \ulow
\eeq
where $\uhigh \in H^2(B_R)$ with
\beq\label{eq:uhigh}
\N{\uhigh}_{H^2_k(B_R)} 
\leq \CHt \N{f}_{L^2(B_R)}
 \quad\text{  for all $k\geq k_0$}
\eeq
and $\ulow \in C^\infty(B_R)$ with
\beq\label{eq:ulow}
\N{(k^{-1}\partial)^\alpha \ulow}_{L^2(B_R)} \leq \Csol(k,2R)\, \big(\CA \big)^{|\alpha|} \N{f}_{L^2(B_R)}\text{ for all $\alpha$ and for all $k\geq k_0$}.
\eeq
\end{theorem}

\paragraph{Understanding the properties of $\uhigh$ and $\ulow$ in Theorem \ref{thm:splitting}.}

Recall that the solution $u$ itself satisfies the bound \eqref{eq:H2}; i.e., $\N{u}_{H^2_k(B_R)} \leq  C kR \N{f}_{L^2(B_R)}$.
Therefore,
\bit
\item[(i)] the bound on $\uhigh$ \eqref{eq:uhigh} is one power of $k$ better than the corresponding bound on $u$ \eqref{eq:H2}, and 
\item[(ii)] the bound on $\ulow$ has the same $k$ dependence as the bound on $u$ \eqref{eq:H2} -- both are governed by $\Csol$ -- although $\ulow$ is $C^\infty$ (with each derivative incurring a power of $k$), and indeed analytic by 
Lemma \ref{lem:analytic}.
\eit
We discuss in \S\ref{sec:why} why both of these points are crucial in proving that the $hp$-FEM does not suffer from the pollution effect. 

%\paragraph{What are $\uhigh$ and $\ulow$?}
%We see in the proof of Theorem \ref{thm:splitting} in \S\ref{sec:proof} that $\uhigh$ corresponds to components of $u$ with frequencies $\geq \lambda k$ and $\ulow$ corresponds to components of $u$ with frequencies $\leq \lambda k$, where $\lambda>1$, and the notion of ``frequencies'' is understood via the Fourier transform.
\paragraph{What are $\uhigh$ and $\ulow$, and why do they satisfy the properties (i) and (ii)?}
We see in the proof of Theorem \ref{thm:splitting} in \S\ref{sec:proof} that $\uhigh$ corresponds to components of $u$ with frequencies $\geq \lambda k$ and $\ulow$ corresponds to components of $u$ with frequencies $\leq \lambda k$, where $\lambda>1$, and the notion of ``frequencies'' is understood via the Fourier transform.
We see in \S\ref{sec:proof_high} below that property (i) above holds because the Helmholtz operator is ``well behaved'' (in a sense made precise below) on frequencies $\geq \lambda k$ with $\lambda>1$. We see in \S\ref{sec:proof_low} that property (ii) above holds because a function with a compactly-supported Fourier transform is analytic.

\subsection{The $hp$-FEM does not suffer from the pollution effect}\label{sec:final}

\begin{theorem}\mythmname{Quasioptimality of the $hp$-FEM}\label{thm:main}
Suppose that $\{\cH_{h,p}\}_{h>0,p\in \mathbb{Z}^+}$ satisfy Assumptions \ref{ass:FEM1} and \ref{ass:FEM2}.
Given $k_0,R_0>0$, there exist $C_1,C_2>0$ (independent of $k$, $R$, $h$, and $p$)
such that 
the following holds.
If $u$ is the solution of the variational problem \eqref{eq:EDPvar2}, $k\geq k_0$,
\beq\label{eq:threshold}
\frac{hk}{p}\leq C_1, \quad\tand\quad p\geq C_2 \log (kR),
\eeq
then the Galerkin solution exists, is unique, and satisfies the quasi-optimal error bound \eqref{eq:qo}. 
\end{theorem}

The pollution effect occurs when no choice of the number of degrees of freedom growing like $(kR)^d$ ensures that the quasi-optimal error bound \eqref{eq:qo_gen} holds with $\Cqo$ independent of $k$ (see \cite[Definition 2.1]{BaSa:00} or \cite[Equation 1.6]{GaSp:22} for more-precise statements of this).
Since the number of degrees of freedom of $\cH_N$ is proportional to $(p/(h/R))^d$, if $h$ and $p$ are chosen so that the inequalities \eqref{eq:threshold} hold with equality, then the number of degrees of freedom of $\cH_N$ is proportional to $(kR)^d$; i.e. Theorem \ref{thm:main} shows that \emph{the $hp$-FEM does not suffer from the pollution effect}.

\bpf[Proof of Theorem \ref{thm:main}]
The plan is to show that there exist $C_1,C_2>0$ such that if $h$, $k$, and $p$ satisfy \eqref{eq:threshold}, then the inequality \eqref{eq:etasuff} holds; the result then follows from Lemma \ref{lem:Schatz}.
By Lemma \ref{lem:Helmholtz_adjoint}, we can consider $\cS^*f$ to be the solution of \eqref{eq:Helmholtz}-\eqref{eq:src};
we then use Theorem \ref{thm:splitting} to split $u$ into $\uhigh$ and $\ulow$, approximate $\uhigh$ using Assumption \ref{ass:FEM1} (with $s=2$), and approximate $\ulow$ by Assumption \ref{ass:FEM2}. By the bounds \eqref{eq:assFEM1} and \eqref{eq:uhigh}, there exists 
$v_{h,p}^{(1)} \in \cH_N$ such that 
\begin{align}
\frac{\|\uhigh - v_{h,p}^{(1)}\|_{H^1_k(B_R)} }
{\N{f}_{L^2(B_R)}}
&\leq \CFEM_1 \left(1 + \frac{hk}{p}\right) \left( \frac{hk}{p}\right) \frac{\|u\|_{H^2_k(B_R)}}{\N{f}_{L^2(B_R)}}\leq \CFEM_1 \left(1 + \frac{hk}{p}\right) \left( \frac{hk}{p}\right) \CHt.
\label{eq:uhighapprox}
\end{align}
The bound \eqref{eq:ulow} implies that $\ulow$ satisfies the conditions of Assumption \ref{ass:FEM2} with 
 $C_2:=\CA$ and $C_1:=
\Csol(k,2R)\N{f}_{L^2(B_R)}$. Therefore, by \eqref{eq:assFEM2} and the bound \eqref{eq:Morawetz_bound} on $\Csol$, there exists $C>0$ and $v_{h,p}^{(2)} \in V_N$ such that 
\beq\label{eq:ulowapprox}
\frac{\|\ulow - v_{h,p}^{(2)}\|_{H^1_k(B_R)} 
}
{\N{f}_{L^2(B_R)}}
\leq 
 \CFEM_2 C
 \left[
\left(\frac{h}{h+\sigma}\right)^p
 + kR\left(\frac{hk}{\sigma p}\right)^p 
\right].
\eeq
Let $v_N:= v_{h,p}^{(1)}+ v_{h,p}^{(2)}$. 
Using the triangle inequality and the decomposition $u=\uhigh+ \ulow$ on $B_R$,
we obtain that 
\beq\label{eq:AHB}
\eta(\cH_{h,p}) \leq  \CFEM_1 \left(1 + \frac{hk}{p}\right) \left( \frac{hk}{p}\right) \CHt + 
 \CFEM_2 C
 \left[
\left(\frac{h}{h+\sigma}\right)^p
 + kR\left(\frac{hk}{\sigma p}\right)^p 
\right].
\eeq
Therefore, to prove the bound \eqref{eq:etasuff}  on $\eta(\cH_{h,p})$ 
it is sufficient to prove that  the right-hand sides of  \eqref{eq:uhighapprox} and \eqref{eq:ulowapprox} are each $\leq \Ccont/4$. To do this, first recall from Lemma \ref{lem:propa} that $\Ccont$ is independent of $k$.
We then choose $C_1$ sufficiently small so that  $C_1\leq\min\{ \widetilde{C}, \sigma\}$ (where $\widetilde{C}$ and $\sigma$ are as in Assumption \ref{ass:FEM2}) and 
\beqs
\CFEM_1 C_1(1+C_1) \CHt\leq \Ccont/4;
\eeqs
observe that, since $hk/p\leq C_1$, this last inequality is the desired bound on the right-hand side of \eqref{eq:uhighapprox}. 
Next let 
\beqs
\theta_1: = \frac{h}{\sigma+ h} \quad\tand\quad \theta_2:= \frac{C_1}{\sigma};
\eeqs
observe that $\theta_1<1$ by definition, and $\theta_2<1$ by the definition of $C_1$.
The right-hand side of  \eqref{eq:ulowapprox} is then bounded by
\beqs
 \CFEM_2 C
\Big[
(\theta_1)^p
 + kR(\theta_2)^p
\Big].
\eeqs
Since $\theta_1,\theta_2<1$, if $p\geq C_2 \log (kR)$ for $C_2$ sufficiently large, then the decay of $(\theta_2)^p$ beats the growth of $kR$; thus with $C_2$ sufficiently large the right-hand side of \eqref{eq:ulowapprox} 
can be made $\leq \Ccont/4$ and the proof is complete.
\epf

\subsection{Discussion of the insight the splitting gives into the pollution effect}
\label{sec:why}

This subsection discusses three natural questions:
\ben
\item Why is the splitting of Theorem \ref{thm:splitting} needed? That is, why does the no-pollution result not just follow from using the bound \eqref{eq:H2} on $u$ itself to bound $\eta(\cH_{h,p})$?
\item How are the properties of $\uhigh$ and $\ulow$ used to prove Theorem \ref{thm:main} (the no-pollution result)?
\item Why does one need $p\to \infty$ to remove the pollution effect?
\een
Regarding 1:~inputting the bound \eqref{eq:H2} on $u$ into the approximation result \eqref{eq:assFEM1} in  Assumption \ref{ass:FEM1}, we obtain that 
\beq\label{eq:etaold}
\eta(\cH_{h,p}) \leq \CFEM_1 \left(1 + \frac{hk}{p}\right) \left( \frac{hk}{p}\right) C kR,
\eeq
which leads to the condition ``$hk^2 R/p$ 
 sufficiently small'' for quasioptimality. The condition ``$hk^2 R$ sufficiently small'' is indeed the observed sharp condition for quasioptimality of the $h$-FEM when $p=1$ -- see, e.g., \cite[Figure 8]{IhBa:95a} -- but in this case the total number of degrees of freedom grows like $(kR)^{2d}$; i.e., the $h$-FEM with $p=1$  suffers from the pollution effect.

Regarding 2 and 3:~$\uhigh$ satisfying the bound \eqref{eq:uhigh}, which is $kR$ better than the corresponding bound \eqref{eq:H2} on $u$, allows us to obtain the condition ``$hk/p$ sufficiently small'' for making the first term on the right-hand side of \eqref{eq:AHB} small, as opposed to the condition ``$hk^2 R/p$ sufficiently small'' for making the right-hand side of \eqref{eq:etaold} small.

When $p$ is fixed, 
to make the second term on the right-hand side of \eqref{eq:AHB} small, we need ``$kR (hk)^p$ sufficiently small'', and the same condition is obtained if we use the approximation result of Assumption \ref{ass:FEM1} to approximate $\ulow$ instead of that of Assumption \ref{ass:FEM2} (i.e., if we ignore the fact that $\ulow$ is analytic, and just use that $\ulow \in H^{s}(B_R)$ for every $s>0$).
The recent polynomial-approximation results of \cite{G1} show that, for fixed $p$,  $\eta(\cH_{h,p})\geq C kR (hk)^p$ (or, more generally, 
$\eta(\cH_{h,p})\geq C \Csol (hk)^p$, where $\Csol$ is defined by \eqref{eq:Csol}). That is, for fixed $p$, the condition ``$kR (hk)^p$ sufficiently small" is the sharp condition for ensuring that $\eta(\cH_{h,p})$ is sufficiently small, and this condition is observed empirically to be the sharp condition required for the Galerkin method to be quasi-optimal with constant independent of $k$ (see, e.g., \cite[Figures 3, 5, and 8]{ChNi:20} for $p=1, 2, 3, 4$); i.e., the $h$-FEM suffers from the pollution effect.

Since $\ulow$ is analytic, we can take $p\to \infty$ in the second term on the right-hand side of \eqref{eq:AHB}, with, importantly, the approximation result of Assumption \ref{ass:FEM2} controlling the dependence of the constant on $p$ (as highlighted in the ``Motivation for Assumption \ref{ass:FEM2}'' paragraph in Section \ref{sec:hpFEM}).
The growth of $kR$ is then removed by the exponential decrease of $(hk/\sigma p)^p$ when $hk/(\sigma p)<1$.
Note that if $kR$ were replaced by $(kR)^M$ for any fixed $M>0$, then this growth would also be 
removed by the exponential decrease of $(hk/\sigma p)^p$ when $hk/(\sigma p)<1$.
Since $kR$ in \eqref{eq:AHB} comes from $\Csol$, the $hp$-FEM results in 
 \cite{MeSa:11, EsMe:12, MePaSa:13, LSW3, LSW4, GLSW1, BeChMe:22}
all involve the assumption that $\Csol$ is polynomially bounded in $kR$.

\section{Recap of results about the Fourier transform and Fourier multipliers}\label{sec:Fourier}

\subsection{The Fourier transform $\mathcal F_k$}\label{sec:Fourier1}

\paragraph{Definition.}
Given $k>0$, let
\beq\label{eq:FT}
\mathcal F_k\phi(\xi) := \int_{\mathbb R^d} \exp\big( -\ri k x \cdot \xi\big)
\phi(x) \, \rd x;
\eeq
i.e., $\mathcal F_k$ is the standard Fourier transform with frequency variable scaled by $k$. The reason for including this scaling is that $\cF_k$ is then tailor-made to work in the weighted Sobolev spaces $H^s_k$ (see \S\ref{sec:Sobolev}), 
 with these spaces the natural spaces to study solutions of the Helmholtz equation.

\paragraph{The Schwartz space and its dual.}
Let $\Schwartz(\Rea^d)$ be the Schwartz space of rapidly decreasing, $C^\infty$ functions; i.e.,
\beqs%\label{eq:Schwartz}
\Schwartz(\Rea^d):= \Big\{ \phi \in C^\infty(\Rea^d) : \sup_{x\in \Rea^d}\big|x^\alpha \partial^\beta \phi(x)\big|<\infty \text{ for all multiindices $\alpha$ and $\beta$}\Big\}.
\eeqs
Let $\Schwartzdual(\Rea^d)$ be the space of continuous linear functionals on $\Schwartz(\Rea^d)$. Recall that $\mathcal F_k: \Schwartz(\Rea^d)\to \Schwartz(\Rea^d)$ 
  (see, e.g., \cite[Page 72]{Mc:00}, \cite[Proposition 13.15]{SaBrHa:19}); then $\mathcal F_k: \Schwartzdual(\Rea^d)\to \Schwartzdual(\Rea^d)$ via the definition that $\langle \mathcal F_k \phi,\psi\rangle:= \langle\phi, \cF_k \psi\rangle$ for $\phi \in \Schwartzdual(\Rea^d)$ and $\psi\in \Schwartz(\Rea^d)$, where $\langle\cdot,\cdot\rangle$ is the duality pairing between $\Schwartzdual(\Rea^d)$ and $\Schwartz(\Rea^d)$. 
In the next subsection we consider $\langle \xi \rangle^s \cF_k \phi$ for $\phi \in \Schwartzdual(\Rea^d)$ and $\langle\xi\rangle := (1+|\xi|^2)^{1/2}$; this is defined as an element of $\Schwartzdual(\Rea^d)$ by $\langle \langle \xi \rangle^s\mathcal F_k \phi,\psi\rangle:= \langle\cF_k \phi, \langle \xi \rangle^s \psi\rangle$ for  $\psi\in \Schwartz(\Rea^d)$; see, e.g., \cite[Propositions 13.14 and 13.17]{SaBrHa:19}.

\paragraph{The properties of $\cF_k$ used in this paper.}
We use the Fourier inversion theorem
\beqs%\label{eq:SCFTinverse}
\mathcal F^{-1}_k\psi(x) := 
\left(\frac{k}{2\pi}\right)^d
\int_{\mathbb R^d} \exp\big( \ri k x \cdot \xi\big)
 \psi(\xi)\, \rd \xi,
\eeqs
the property
\beq\label{eq:FTderivative}
\mathcal{F}_k \big( \big( -\ri \hsc \partial\big)^\alpha \phi\big)(\xi) = \xi^\alpha \,\mathcal{F}_k \phi(\xi),
\eeq
and Plancherel's theorem
\beq\label{eq:Plancherel}
\N{\phi}_{L^2(\Rea^d)} = \left(\frac{k}{2\pi}\right)^{d/2}\N{\mathcal{F}_k \phi}_{L^2(\Rea^d)}.
\eeq

\subsection{Sobolev spaces weighted with $k$}\label{sec:Sobolev}

The natural spaces in which to study solutions of the Helmholtz equation are Sobolev spaces 
with derivatives weighted with $k$, as in \eqref{eq:weighted_norms}.
On $\Rea^d$ these are naturally defined using $\mathcal{F}_k$.
For $s\in\Rea$, let
\beqs%\label{eq:Sobolev}
H_k ^ s (\Rea^d):= \Big\{ u\in \Schwartzdual(\mathbb R^d), \; \langle \xi \rangle^s 
\mathcal F_k u \in  L^2(\mathbb R^d) \Big\}, \quad\text{ where }\langle \xi \rangle := (1+|\xi|^2)^{1/2}, 
\eeqs
and let
\beq\label{eq:Hhnorm}
\vertiii{u}_{H_k^s(\Rea^d)} ^2 := 
\left(\frac{k}{2\pi}\right)^d
 \int_{\Rea^d} \langle \xi \rangle^{2s}
 |\mathcal F_k u(\xi)|^2 \, \rd \xi.
\eeq
Because of \eqref{eq:FTderivative}, up to dimension-dependent constants, $\vertiii{u}_{H_k^s(\Rea^d)}$ defined by \eqref{eq:Hhnorm} is equivalent to $\| u \|_{H^s_k(\Rea^d)}$ defined by \eqref{eq:weighted_norms} (with $B_R$ replaced by $\Rea^d$).
To be more precise about the norm equivalence, if $C_j=C_j(s,d)>0$, $j=1,2,$ are such that 
\beq
C_1 \sum_{|\alpha|\leq s} \xi^{2\alpha} \leq (1+ |\xi|^2)^{s} \leq C_2 \sum_{|\alpha|\leq s} \xi^{2\alpha},
\label{eq:normequiv}
\,\,\text{ then }\,\,
\sqrt{C_1} \N{u}_{H^{s}_k(\Rea^d)} \leq \vertiii{u}_{H^s_k(\Rea^d)}\leq  \sqrt{C_2} \N{u}_{H^s_k(\Rea^d)}.
\eeq

\subsection{Fourier multipliers}\label{sec:Fouriersymbol}

\paragraph{The basic idea.} The \emph{Fourier multiplier} given by a function $a$ is
\beq \label{eq:OphF}
\big(a(\hsc D)v\big)(x):= \mathcal{F}^{-1}_k \big(a(\cdot)(\mathcal{F}_k v)(\cdot)\big)(x);
\eeq
i.e., we multiply the Fourier transform of $v$ by $a$, and then apply the inverse Fourier transform.
The rationale for the notation $a(\hsc D)$ is that, by \eqref{eq:FTderivative}, if $D:= -\ri \partial$ then $\mathcal{F}_k$ maps $\hsc D$ to $\xi$.

Our motivation for studying Fourier multipliers is that (i) the operator $-k^{-2} \Delta -1$ is one -- by the derivative rule \eqref{eq:FTderivative}, $-k^{-2}\Delta-1 = p(k^{-1}D)$ where $p(\xi):=|\xi|^2-1$, and (ii) the functions $\uhigh$ and $\ulow$ in our proof of Theorem \ref{thm:splitting} are defined by Fourier multipliers acting on $u$ (see \eqref{eq:uhighlow} below).

\paragraph{A natural class of functions for which Fourier multipliers are well-defined on Sobolev spaces.}
We say that $a$ is a \emph{Fourier symbol of order $m$} if 
there exists $C>0$ such that 
\beq\label{eq:FSm}
|a(\xi) | \leq C \langle \xi \rangle^{m} 
\quad\tfa \xi \in \Rea^d, 
\eeq
where recall that $\langle \xi \rangle:=(1+|\xi|^2)^{1/2}$.
We use the (non-standard) notation that $a \in (FS)^m$.

\begin{example}\mythmname{Examples of Fourier symbols and multipliers}\label{ex:Fouriermult}

(i) If $a(\xi)=1$, then $a\in (FS)^0$ with $a(\hsc D)v(x)= v(x)$ 
(since $\cF^{-1}_k\cF_k=I$).

(ii) If $p(\xi):= |\xi|^2 -1$ then $p\in (FS)^2$ with $p(\hsc D)v = (- k^{-2} \Delta -1)v$.

(iii) If $\chi$ is bounded and has compact support, then $\chi\in (FS)^{-N}$ for all $N\geq 1$ and $1-\chi \in (FS)^0$.
\end{example}

\begin{theorem}\mythmname{Composition and mapping properties of Fourier multipliers}
\label{thm:basicPFourier} 
If $a \in (FS)^{m_a}$ and $b \in (FS)^{m_b}$ then the following properties hold.

(i) $ab\in (FS)^{m_a+m_b}$.

(ii) $a(\hsc D) b(\hsc D)= (ab)(\hsc D)=b(\hsc D) a(\hsc D)$.

(iii) $a(\hsc D): H_k^{s}(\Rea^d)\to H^{s-m_a}_k(\Rea^d)$ and, with $C$ the constant in \eqref{eq:FSm}, for all $s\in \Rea$ and $k >0$,
\beq\label{eq:mapping}
\vertiii{a(\hsc D)}_{H^s_k(\Rea^d) \to H^{s-m_a}_k(\Rea^d)} \leq C;
\eeq
i.e., $a(\hsc D)$ is bounded uniformly in both $k$ and $s$ as an operator from $H_k^s$ to $H_k^{s-m_a}$.
\end{theorem}

\bpf
(i) This follows directly from the definition \eqref{eq:FSm}. 

(ii) By the definition \eqref{eq:OphF},
\begin{align*}
a(\hsc D) \big(b(\hsc D)v\big)(x) = \mathcal{F}_k^{-1}\big( a(\cdot)(\mathcal{F}_k(b(\hsc D)v))(\cdot)\big)(x)
&= \mathcal{F}_k^{-1}\big( a(\cdot) b(\cdot) (\mathcal{F}_k v)(\cdot) \big)(x) \\
&\hspace{-2cm}= (ab)(\hsc D)v(x)\\
&\hspace{-2cm}= (ba)(\hsc D)v(x) = 
b(\hsc D) \big( a(\hsc D)\big)v(x).
\end{align*}

(iii) 
By the definitions of $\vertiii{\cdot}_{H^s_k(\Rea^d)}$ \eqref{eq:Hhnorm} and $a(\hsc D)v$ \eqref{eq:OphF},
\begin{align*}
\vertiii{a(\hsc D)v}^2_{H^{s-m_A}_k(\Rea^d)} 
&=  \left(\frac{k}{2\pi}\right)^d\int_{\Rea^d} \langle \xi\rangle^{2(s-m_A)} \big|a(\xi)\big(\mathcal{F}_k v\big)(\xi)\big|^2 \, \rd \xi\\
&\leq C  \left(\frac{k}{2\pi}\right)^d \int_{\Rea^d} \langle \xi\rangle^{2(s-m_A)}\langle \xi\rangle^{2m_A} \big|\big(\mathcal{F}_k v\big)(\xi)\big|^2 \, \rd \xi
= C\vertiii{v}_{H^s_k(\Rea^d)}^2.
\end{align*}
\epf

The key result from this section used in the proof of the bound \eqref{eq:uhigh} on $\uhigh$ is the following.

\begin{theorem}\mythmname{Factoring an ``elliptic'' Fourier multiplier out of another}\label{thm:elliptic}
Suppose that $a\in (FS)^{m_a}$, $b \in (FS)^{m_b}$, and there exists $c>0$ such that 
\beq\label{eq:b_elliptic}
|b(\xi)|\geq c\langle\xi\rangle^{m_b} \quad\tfor \xi \in\supp\, a.
\eeq
Then
\beq\label{eq:elliptic1}
a(\hsc D) = q(\hsc D) b(\hsc D)
\eeq
where $q\in (FS)^{m_a-m_b}$ is defined by 
\beq\label{eq:q}
q(\xi) := \frac{a(\xi)}{b(\xi)}.
\eeq
\end{theorem}

\bpf
The fact that $q$ defined by \eqref{eq:q} is in $(FS)^{m_a-m_b}$ follows from the fact that $a\in (FS)^{m_a}$ and the bound \eqref{eq:b_elliptic}. 
The result \eqref{eq:elliptic1} then follows from Part (ii) of Theorem \ref{thm:basicPFourier}.
\epf

We now combine Theorem \ref{thm:elliptic} and the mapping property \eqref{eq:mapping} to obtain the following result (which we use in the proof of Corollary \ref{cor:H2}); we highlight that a similar combination of these results in used 
 in \S\ref{sec:proof} in the proof of the bound \eqref{eq:uhigh} on $\uhigh$.
 
\begin{corollary}\mythmname{Elliptic regularity}\label{cor:elliptic_regularity}
If $(-k^{-2}\Delta + 1) v \in L^2(\Rea^d)$ then $v\in H^2(\Rea^d)$ with
\beqs%\label{eq:elliptic_reg}
\vertiii{v}_{H^2_k(\Rea^d)}\leq \N{(-k^{-2}\Delta + 1 )v}_{L^2(\Rea^d)}.
\eeqs
\end{corollary}

\bpf
We apply Theorem \ref{thm:elliptic} with $a(\xi)=1$ and $b(\xi)= |\xi|^2+1$, so that $a(k^{-1}D)= I$ and $b(k^{-1}D)= -k^{-2}\Delta +1$. The theorem implies that $q(\xi):=\langle \xi\rangle^{-2}\in (FS)^{-2}$ (this is also clear from the definition \eqref{eq:FSm}). Then,
by the mapping property \eqref{eq:mapping},
\beqs
\vertiii{v}_{H^2_k(\Rea^d)} = \vertiii{q(k^{-1}D) (-k^{-2}\Delta +1)v }_{H^2_k(\Rea^d)} \leq \N{(-k^{-2}\Delta +1)v}_{L^2(\Rea^d)}.
\eeqs
\epf

\section{Proof of Theorem \ref{thm:splitting} using only the material in \S\ref{sec:Fourier}}\label{sec:proof}

\subsection{Definition of high- and low- frequency cut-offs}

Let 
\beq\label{eq:chi}
\chi_\lambda(\xi) = 1_{|\xi|\leq \lambda}(\xi) = 
\begin{cases}
1 & \text{ for } |\xi|\leq \lambda,\\
0 & \text{ for } |\xi|>\lambda.
\end{cases}
\eeq
We define the low-frequency cut-off $\Pilow$ by
\beq
\Pilow := \chi_\lambda(k^{-1}D); \,\,\quad\text{i.e.,} \,\,\quad
\label{eq:Pilow}
\Pilow v = \cF_k^{-1}\big( \chi_\lambda(\cdot) (\cF_k v)(\cdot)\big),
\eeq
by the definition \eqref{eq:OphF} of a Fourier multiplier.
We see that $\Pilow$ acting on a function $v$ returns the frequencies of $v$ that are $\leq \lambda$, hence why we call it a low-frequency cut-off. 
\footnote{The definition of $\cF_k$ \eqref{eq:FT} implies that if $\cF_k w$ is supported on $|\xi|\leq \lambda$, then the ``standard'' Fourier transform (i.e., with the transform variable not scaled by $k$) of $w$ is supported for $|\zeta|\leq \lambda k$.} % (where $\zeta$ is the ``standard'' Fourier variable).}
We define the high-frequency cut-off $\Pihigh$ by 
\beqs%\label{eq:Pihigh}
\Pihigh: = I- \Pilow = (1-\chi_\lambda)(k^{-1}D);
\eeqs
i.e., $\Pihigh$ acting on a function $v$ returns the frequencies of $v$ that are $\geq \lambda$.

\subsection{Definition of $\uhigh$ and $\ulow$ via the frequency cut-offs}\label{sec:defuhighulow}

Let $u$ be as in Theorem \ref{thm:splitting}; i.e., $u$ is the outgoing solution of the Helmholtz equation \eqref{eq:Helmholtz}. Let  $\varphi \in C^\infty_{\rm comp}(\Rea^d, [0,1])$ be equal to one on $B_{R}$ and vanish outside $B_{2R}$, and set
\beq\label{eq:uhighlow}
\ulow := \big(\Pilow (\varphi u)\big)\big|_{B_R}
\quad\tand\quad
\uhigh := \big(\Pihigh(\varphi u)\big)\big|_{B_R}.
\eeq
Since $\Pilow+ \Pihigh=I$, these definitions imply that, on $B_R$, $ \ulow + \uhigh =\varphi u= u$; i.e., \eqref{eq:splitting} holds.
This splitting contains the arbitrary parameter $\lambda$; we fix this when proving the bound \eqref{eq:uhigh} on $\uhigh$.

\subsection{Proof of the bound \eqref{eq:ulow} on $\ulow$}\label{sec:proof_low}

The idea of the proof, in short, is that \emph{a function with a compactly-supported Fourier transform is analytic}.  
Plancherel's theorem \eqref{eq:Plancherel} and the derivative property \eqref{eq:FTderivative} 
and the definition of $\Pilow$ \eqref{eq:Pilow}
imply that
\begin{align}%\nonumber
  \N{(k^{-1}\pa)^\alpha\big( \Pilow \varphi u\big)}_{L^2(\Rea^d)}=
  \left(\frac{k}{2\pi}\right)^{d/2}
  \N{(\cdot)^\alpha \mathcal{F}_k\big( \Pilow \varphi u\big)(\cdot)}_{L^2(\Rea^d)} 
=    \left(\frac{k}{2\pi}\right)^{d/2}\N{(\cdot)^\alpha \chi_\lambda(\cdot)\mathcal{F}_k(\varphi u)(\cdot)}_{L^2(\Rea^d)}. \label{eq:VPS1}
  \end{align}
The definition of $\chi_\lambda$ \eqref{eq:chi} implies that 
\beqs
\big| \xi^\alpha \chi_\lambda(\xi)\big| \leq \lambda^{|\alpha|} \quad \tfa \xi \in \Rea^d.
\eeqs  
Using this fact in \eqref{eq:VPS1} and then (in this order) Plancherel's theorem \eqref{eq:Plancherel}, the fact that $\varphi=0$ outside $B_{2R}$ and $\varphi \leq 1$ inside $B_{2R}$, and the definition of $\Csol$ \eqref{eq:Csol}, we find that
\begin{align}\nonumber
  \N{(k^{-1}\pa)^\alpha\big( \Pilow \varphi u\big)}_{L^2(\Rea^d)}
  \leq 
    \left(\frac{k}{2\pi}\right)^{d/2}
   \lambda^{|\alpha|} \N{\mathcal{F}_k(\varphi u)}_{L^2(\Rea^d)}%\\ \nonumber
    &\leq    \lambda^{|\alpha|}\N{\varphi u}_{L^2(\Rea^d)}\\ %\nonumber
        &\hspace{-2.5cm}\leq    \lambda^{|\alpha|}\N{u}_{L^2(B_{2R})}\leq    \lambda^{|\alpha|}\Csol(k, 2R) \N{f}_{L^2(B_R)}.\label{eq:VPS2}
  \end{align}
By the definition \eqref{eq:uhighlow} of $\ulow$,
\beqs
\N{(k^{-1}\pa)^\alpha \ulow}_{L^2(B_R)} = \N{(k^{-1}\pa)^\alpha \big(\Pilow \varphi u\big)}_{L^2(B_R)} \leq  \N{(k^{-1}\pa)^\alpha \big(\Pilow \varphi u\big)}_{L^2(\Rea^d)};
\eeqs
the bound \eqref{eq:ulow} then follows from \eqref{eq:VPS2} with  $\CA:= \lambda$. 
We make two remarks.
\bit
\item
The fact that $\Pilow \varphi u\in H^s_k(\Rea^d)$ for all $s$ follows from the fact that $\chi_\lambda \in (FS)^{-N}$ for every $N>0$ (from Part (iii) of Example \ref{ex:Fouriermult}) and the mapping property in Part (iii) of Theorem \ref{thm:basicPFourier}; we give the direct proof above, however, to have explicit control on the constants in the bound (to show that $\Pilow \varphi u$ is actually analytic).
\item
The fact that $\ulow$ comes from a function with a compactly-supported Fourier transform, and hence automatically is analytic, is one of the advantages of the current splitting compared to the original splitting in \cite{MeSa:10}; see the discussion in Section \ref{sec:MS}.
\eit
%We note that the fact that $\Pilow \varphi u\in H^s_k(\Rea^d)$ for all $s$ follows from the fact that $\chi_\lambda \in (FS)^{-N}$ for every $N>0$ (from Part (iii) of Example \ref{ex:Fouriermult}) and the mapping property in Part (iii) of Theorem \ref{thm:basicPFourier}; we give the direct proof above, however, to have explicit control on the constants in the bound (to show that $\Pilow \varphi u$ is actually analytic).

\subsection{Proof of the bound \eqref{eq:uhigh} on $\uhigh$}\label{sec:proof_high}

The idea of the proof is to use Theorem \ref{thm:elliptic}, using the fact that \emph{the Helmholtz operator is an elliptic Fourier multiplier (in the sense of \eqref{eq:b_elliptic}) on high frequencies}, and thus, if $\lambda$ is sufficiently large, on the support of the high-frequency cut-off $\Pihigh$.

By the definition of $\uhigh$ \eqref{eq:uhighlow} and the equivalence  of $\|\cdot\|_{H^2_k}$ and $\vertiii{\cdot}_{H^2_k}$ described in \eqref{eq:normequiv},
\beqs
\N{\uhigh}_{H^2_k(B_R)} = \N{\Pihigh(\varphi u)}_{H^2_k(B_R)} \leq  \N{\Pihigh(\varphi u)}_{H^2_k(\Rea^d)} \lesssim 
\vertiii{\Pihigh(\varphi u)}_{H^2_k(\Rea^d)},
\eeqs
where we use the notation that $A\lesssim B$ if there exists $C>0$, independent of $k$ and $R$, such that $A \leq CB$.
It is therefore sufficient to prove that
\beq\label{eq:uhighFourierbound}
\vertiii{\Pihigh(\varphi u)}_{H^2_k(\Rea^d)}
\lesssim \N{f}_{L^2(B_R)} \quad \tfa k\geq k_0.
\eeq

\ble\label{lem:Fourierproof}\mythmname{``Ellipticity'' of $p(\xi):=|\xi|^2-1$ when $|\xi|\geq\lambda>1$}
If $\lambda \geq \lambda_0>1$ then there exists $C>0$ such that 
\beq\label{eq:VPS3}
\big||\xi|^2-1\big|\geq C \langle\xi\rangle^2 \quad \tfor |\xi|\geq  \lambda.
\eeq
\ele

\bpf
It is straightforward to check that \eqref{eq:VPS3} holds with $C:= (1 + 2(\lambda_0^2-1)^{-1})^{-1}$.
\epf

\begin{corollary}\mythmname{Factoring out $p(k^{-1}D)$ from $(1-\chi_\lambda)(k^{-1}D)$}\label{cor:key}
Let $p(\xi):=|\xi|^2-1$. 
If $\lambda \geq \lambda_0>1$, then there exists $q\in (FS)^{-2}$ such that 
\beqs
(1-\chi_\lambda)(k^{-1}D) = q(k^{-1}D) \, p(k^{-1}D).
\eeqs
\end{corollary}

\bpf
By the definition \eqref{eq:chi} of $\chi_\lambda$,
$\supp (1-\chi_\lambda)= \{\xi : |\xi|\geq \lambda\}$. 
The result follows by applying Theorem \ref{thm:elliptic} with $a(\xi):= (1-\chi_\lambda)(\xi)$ and $b(\xi):= p(\xi)=|\xi|^2-1$, since 
\eqref{eq:VPS3}
 implies that the inequality \eqref{eq:b_elliptic} holds (i.e., $b$ is elliptic on $\supp \,a$).
\epf

We now use Corollary \ref{cor:key} and the mapping property \eqref{eq:mapping} to prove the bound \eqref{eq:uhighFourierbound}. 
By the definition of $\uhigh$ \eqref{eq:uhighlow} and Corollary \ref{cor:key},
\begin{align}\nonumber 
\vertiii{\Pihigh (\varphi u)}_{H^2_k(\Rea^d)} &= \vertiii{\big(1- \chi_\mu\big)(\hsc D)(\varphi u)}_{H^2_k(\Rea^d)}
= \vertiii{q(\hsc D) \,p(\hsc D)(\varphi u)}_{H^2_k(\Rea^d)}%\label{eq:Fourierproof1}
\end{align}
with $q\in (FS)^{-2}$, and, by the mapping property \eqref{eq:mapping}, 
\beqs
\vertiii{q(\hsc D) \,p(\hsc D)(\varphi u)}_{H^2_k(\Rea^d)}
\lesssim \N{p(\hsc D)(\varphi u)}_{L^2(\Rea^d)};
\eeqs
the key point is that we now have the Helmholtz operator $p(k^{-1}D)$ on the right-hand side, and we can start to relate this side to 
$f= p(k^{-1}D)u$. For brevity, let $P:= -k^{-2}\Delta -1 = p(k^{-1}D)$. Then 
\begin{align}
\vertiii{\Pihigh (\varphi u)}_{H^2_k(\Rea^d)} 
\lesssim \N{P(\varphi u)}_{L^2(\Rea^d)} &= 
\N{\varphi P u + [P, \varphi] u }_{L^2(\Rea^d)}\leq \N{f}_{L^2(\Rea^d)}  + \N{[P, \varphi] u }_{L^2(\Rea^d)},
\label{eq:Fourierproof2}
\end{align}
where we have used the fact that $\varphi \equiv 1$ on $\supp\, f$, and where the commutator $[A,B]$ is defined as $AB-BA$.
By direct calculation,
\beq\label{eq:Fourierproof3}
[P,\varphi]u= -k^{-2}\big( u\Delta \varphi +2 \nabla \varphi \cdot\nabla u),
\quad
\text{ so that }
\quad
\big\| [P,\varphi]u\big\|_{L^2(\Rea^d)} \lesssim (kR)^{-1} \N{u}_{H^1_k(B_{2R})},
\eeq
where we have used that the definition of $\varphi$ in \S\ref{sec:defuhighulow} implies that $|\nabla \varphi| \sim R^{-1}$ and $|\Delta \varphi| \sim R^{-2}$. Therefore, by combining \eqref{eq:Fourierproof2} and \eqref{eq:Fourierproof3}, and using \eqref{eq:H2} (with $R$ replaced by $2R$), 
we have 
\beqs
\N{\Pihigh (\varphi u)}_{H^2_k(\Rea^d)} 
\lesssim   \N{f}_{L^2(\Rea^d)} +(kR)^{-1} \N{u}_{H^1_k(B_{2R})}\lesssim \N{f}_{L^2(\Rea^d)},
\eeqs
which is the result \eqref{eq:uhighFourierbound}.

\subsection{Discussion of the original proof of Theorem \ref{thm:splitting} in \cite{MeSa:10}}\label{sec:MS}

The original proof of Theorem \ref{thm:splitting} in \cite[\S3.2]{MeSa:10} is also based on the idea of frequency cut-offs using the indicator function \eqref{eq:chi}.
 However, in \cite[\S3.2]{MeSa:10} the frequency cut-offs are 
 applied to the data $f$ (as opposed to $\varphi u$ in our case). 
Furthermore, the analysis in \cite[\S3.2]{MeSa:10} is based on writing the solution $u$ to the Helmholtz equation \eqref{eq:Helmholtz} satisfying the Sommerfeld radiation condition \eqref{eq:src} as 
\beq\label{eq:convolution}
u(x) = \int_{B_R}\Phi_k(x,y) \, f(y) \, \rd y,
\eeq
where $\Phi_k(x,y)$ is defined by \eqref{eq:fund}. 
The proof in  \cite[\S3.2]{MeSa:10} then considers the function 
\beqs
v_\mu(x) := \int_{B_R}\Phi_k(x,y) \, \mu(|x-y|)\, f(y) \, \rd y, 
\eeqs
where $\mu \in C^\infty_{\rm comp}(\Rea)$ with $\mu|_{[0,2R]}=1$, $\supp\, \mu \subset[0, 4R]$, and additionally the first and second derivatives of $\mu$ satisfying certain bounds (see \cite[Equation 3.27]{MeSa:10}). This definition of $\mu$ implies that $v_\mu |_{B_R} = u|_{B_R}$. The advantage of studying $v_\mu$ instead of $u$ is that $v$ is the convolution of a compactly-supported kernel with $f$; thus $v$ has compact support and its Fourier transform is well defined. (In contrast, in the present paper, we first multiply $u$ by $\varphi$ to ensure that $\cF_k(\varphi u)$ makes sense.)
The analysis in \cite[\S3.2 and Appendix A]{MeSa:10} then proceeds by studying the compactly-supported integral kernel in $v_\mu$ and the Fourier transform of this kernel, with identities and bounds on Bessel and Hankel functions needed when both $d=2$ and $d=3$.

\subsection{Generalising Theorem \ref{thm:splitting} to more-complicated Helmholtz problems}\label{sec:LSW}

The proof of Theorem \ref{thm:splitting} above can be generalised to more-complicated Helmholtz problems, involving variable coefficients and obstacles, using pseudodifferential operators. Indeed, whereas the operator $-k^{-2} \Delta -1$ is a Fourier multiplier (by Part (ii) of Example \ref{ex:Fouriermult}), 
its variable coefficient analogue $-k^{-2} \nabla \cdot (A \nabla) -n$ is not; i.e., one cannot write the Fourier transform of $-k^{-2} \nabla \cdot (A \nabla v) -n v$ in terms of the Fourier transform of $v$. 
Nevertheless, \cite{LSW3} proves Theorem \ref{thm:splitting}, using the same ideas in the proof above, with the Helmholtz equation \eqref{eq:Helmholtz} replaced by
\beq\label{eq:Helmholtz_variable}
k^{-2}\nabla \cdot (A\nabla u ) + n u = -f \quad\tin \Rea^d.
\eeq
This is achieved by replacing the Fourier multipliers in \S\ref{sec:Fourier} by \emph{pseudodifferential operators} (indeed, recall that 
one of the motivations for the development of the theory of pseudodifferential operators in the 1960s was to study variable-coefficient PDEs, such as \eqref{eq:Helmholtz_variable}, using Fourier analysis).

We therefore hope that the material in \S\ref{sec:Fourier} and \S\ref{sec:proof}, combined with the introductory discussion of pseudodifferential operators in \cite{LSW3}, provides for the interested reader a ``bridge'' 
between Fourier multipliers and pseudodifferential operators. In particular, the pseudodifferential generalisation of the key result of Theorem \ref{thm:elliptic} is the so-called ``elliptic parametrix''.

These ideas can also be used to prove an analogue of Theorem \ref{thm:splitting} when the Helmholtz equation (with or without variable coefficients) is posed not in $\Rea^d$ but outside an impenetrable obstacle. In this case the Fourier transform \eqref{eq:FT} is no longer available; nevertheless the \emph{functional calculus} (essentially the idea of expanding in eigenfunctions of the differential operator) can be used to define Fourier-type transforms, tailored to the particular problem. The ideas of the proof of Theorem \ref{thm:splitting} can then be implemented in this situation (albeit with more technicalities, and the requirement that the obstacle is analytic); see \cite{LSW4, GLSW1}.

\bre[The smoothness of the frequency cut-offs]
The theory of pseudodifferential operators is easiest when the symbols (i.e., the generalisation of the Fourier symbols in \eqref{eq:FSm}) are $C^\infty$. To ensure this, \cite{LSW3} assumes that the coefficients $A$ and $n$ in \eqref{eq:Helmholtz_variable} are $C^\infty$, and uses a smooth frequency cut-off; i.e., the indicator function in \eqref{eq:chi} is replaced by
$\chi_\lambda \in C^\infty_{\rm comp}(\Rea^d;[0,1])$ such that $\chi_\lambda(\xi)=1$ for $|\xi|\leq \lambda$ and $\chi_\lambda(\xi)=0$ for, say, $|\xi|\geq 2\lambda$. The reader can check that the proof in \S\ref{sec:proof} goes through as before with this smooth cut-off (with $\CA$ in \S\ref{sec:proof_low} changed from $\lambda$ to $2\lambda$).
\ere

\appendix

\section{Proof of Theorem \ref{thm:Morawetz}}\label{app:Morawetz}

\ble\mythmname{Morawetz identity for the Helmholtz operator \cite[\S I.2]{Mo:75}}\label{lem:M1}
If
\beq\label{eq:LM}
\opL v:= k^{-2}\Delta v+ v
\quad\tand\quad
\cM_{\beta,\alpha} v:= \bx\cdot \gv - \ri k \beta v + \alpha v,
\eeq
with $\beta$ and $\alpha$ real-valued $C^1$ functions, then
\begin{align}\nonumber 
\hspace{-0.15cm}2 \Re \big(\overline{\cM_{\beta,\alpha} v } \,\cL v \big) 
= &\, \nabla \cdot \Big[ 2 k^{-1}\Re\big(\overline{\cM_{\beta,\alpha} v}\, k^{-1}\nabla v \big) + \big(\nvs - k^{-2}\ngvs\big) \bx\Big] \\
&- 2 \Re\big( \overline{v}\,(\ri  \nabla \beta + k^{-1}\nabla \alpha)\cdot k^{-1}\nabla v \big)
-\big(d-2 \alpha \big) \nvs  -\big(2\alpha -d+2\big)k^{-2}|\gv|^2. \label{eq:morid1}
\end{align}
\ele

\bpf
This follows in a straightforward (but slightly involved) way by expanding the divergence on the right-hand side; for this done step-by-step, see, e.g., \cite[Proof of Lemma 2.1]{SpKaSm:15}.
\epf

The idea of the proof of Theorem \ref{thm:Morawetz} is to integrate the identity \eqref{eq:morid1} over $\Rea^d$ with $v=u$, $\alpha=(d-1)/2$, and $\beta$ defined piecewise as $\beta= R$ for $r\leq R$ and $\beta= r$ for $r\geq R$. The choice $\beta=$ constant and $\alpha=(d-1)/2$ means that the non-divergence terms on the right-hand side of \eqref{eq:morid1} become $-|u|^2 - k^{-2}|\nabla u|^2$; this is where we get $\|u\|_{H^1_k(B_R)}^2$ from. The choice $\beta=r$ deals with the contribution from infinity (although this is not immediately clear from \eqref{eq:morid1}). We therefore first look at the special case of \eqref{eq:morid1} with $\beta=r$.

\ble\mythmname{Special case of \eqref{eq:morid1} \cite[Equation 1.2]{MoLu:68}}
With $\opL v$ and $\cM_{\beta,\alpha} v$ as in Lemma \ref{lem:M1}, show that if $\alpha\in \Rea$, then, with $v_r=\bx\cdot \gv/r$,
\bal\nonumber
2\Re\big( \overline{\cM_{r,\alpha} v} \opL v\big) =  &\,\nabla \cdot \Big[2k^{-1}\Re \big(\overline{\cM_{r, \alpha} v} \,k^{-1}\gv\big)+ \left(\nvs - k^{-2}\ngvs \right)\bx\Big]-\big| k^{-1}v_r -\ri v\big|^2 \\&\quad + \big(2\alpha -(d-1)\big)\big( \nvs - k^{-2}\ngvs\big) - k^{-2}\big(\ngvs -\nvrs\big).
\label{eq:ml2d}
\end{align}
\ele

\bpf
This follows from \eqref{eq:morid1} by choosing $\beta=r$ and writing the term involving $\nabla \beta$ as
\beqs
-2\Re{\left(\ri\overline{v} k^{-1}v_r\right)} = -2\Re{\big(\overline{(-\ri v)} k^{-1}v_r\big)} = \abs{v}^2 + k^{-2}\abs{v_r}^2 - \big|k^{-1}v_r - \ri v\big|^2
\eeqs
using the identity $-2\Re{(z_1\overline{z_2})} =  |z_1|^2 + |z_2|^2-|z_1 + z_2|^2 $.
\epf

To integrate \eqref{eq:ml2d} over $\Rea^d\setminus B_R$, we integrate it over $B_{R_1}\setminus B_R$ and then send $R_1\to \infty$. In preparation for this, we look at the boundary term on $\partial B_{R_1}$.

\ble\label{lem:MLint}
Let
\beq\label{eq:Q}
 \bQ_{r,\alpha}(v):=
2k^{-1}\Re \big(\overline{\cM_{r, \alpha} v} \,k^{-1}\gv\big)+ \left(\nvs - k^{-2}\ngvs \right)\bx.
\eeq
If $u$ is an outgoing solution of $\opL u=0$ in $\Rea^d\setminus \overline{B_{R_0}}$,
 then, for all $\alpha\in \Rea$,
\beq\label{eq:N0}
\int_{\Gamma_{\partial B_{R_1}}} 
\bQ_{R_1,\alpha}(u)
\cdot\hatx
 \tendo \quad \tas R_1\tendi.
\eeq
\ele

\bpf
By the definitions of $\bQ_{r,\alpha}(v)$ \eqref{eq:Q} and $\cM_{r, \alpha}v$ \eqref{eq:LM}, 
\begin{align}\nonumber
Q_{r,\alpha}(u) \cdot \widehat{x} &= rk^{-2}\left(2\Re{\left(u_r\overline{\left(u_r - \ri ku + \frac{\alpha}{r}u\right)}\right)} + k^2|u|^2 -  \abs{\nabla u}^2\right)\\ \nonumber
&= rk^{-2}\left(|u_r|^2 + 2 \Re{\left(u_r\overline{\left(- \ri ku + \frac{\alpha}{r}u\right)}\right)}
+ k^2 |u|^2- \big(|\nabla u|^2 - |u_r|^2\big) 
\right)\\ \label{eq:VPS6}
&=rk^{-2}\left( \frac{|\mathcal{M}_{r,\alpha}u|^2}{r^2} + \alpha^2 \frac{|u|^2}{r^2} - \big(|\nabla u|^2 - |u_r|^2\big)\right),
\end{align}
where we have again used the identity $2\Re{(z_1\overline{z_2})} = |z_1 + z_2|^2 - |z_1|^2 - |z_2|^2$. 
We now claim that the term in large brackets in \eqref{eq:VPS6} is $O(r^{-d-1})$; if this is true, then
\beqs
\int_{\Gamma_{\partial B_{R_1}}} 
\bQ_{R_1,\alpha}(u)
\cdot\hatx= O\left(\frac{1}{R_1}\right) \tas R_1\tendi,
\eeqs
and thus \eqref{eq:N0} follows.
By the Atkinson--Wilcox expansion \eqref{eq:AtWi_intro}/\eqref{eq:AtWi},
$|u|^2 = O(r^{1-d})$ and $r^{-2}|\mathcal{M}_{r,\alpha}u|^2= O(r^{-d-1})$. To prove the result, therefore, we only need to show that $|\nabla u|^2 - |u_r|^2= O(r^{-d-1})$.
The quantity $\ngtus$ equals $|\nabla_S u|^2$ where $\nabla_S$ is the surface gradient on $|x|=r$, 
which satisfies $\nabla_S u = \gu - \widehat{x}u_r$. This differential operator is equal to $1/r$ multiplied by an operator acting only on $\widehat{x}$, i.e. the angular variables; thus $|\nabla_Su|^2$ is $O(r^{-d-1})$ and the proof is complete.
\epf

We now integrate \eqref{eq:ml2d} over $B_{R_1}\setminus B_R$, send $R_1\to \infty$, and obtain an inequality involving the boundary term on $\partial B_R$.

\ble\label{lem:M4}
If $u$ is an outgoing solution of $\cL u=0$ in $\Rea^d\setminus \overline{B_{R_0}}$, for some $R_0>0$, 
then, for $R>R_0$, 
\beq\label{eq:2.1}
\int_{\partial B_R} \bQ_{R,(d-1)/2}(u)
\cdot\hatx \leq 0.
\eeq
\ele

\bpf
We integrate the identity \eqref{eq:ml2d} over $B_{R_1}\setminus B_R$, where $R_1>R$, with $v=u$ and then use the divergence theorem $\int_D \nabla \cdot F = \int_{\partial D} F$. The divergence theorem is valid for $F\in C^\infty(\overline{D})$ and $D$ Lipschitz (see, e.g., \cite[Theorem 3.34]{Mc:00}); we can use it here since, by elliptic regularity, $u\in C^\infty(\overline{B_{R_1}\setminus B_R})$ (see, e.g., \cite[Theorem 4.16]{Mc:00}).
This results in 
\begin{align}\nonumber
&\int_{\partial B_{R_1}} Q_{R_1,\alpha}(u) \cdot \widehat{x} - \int_{\partial B_R} Q_{R,\alpha}(u) \cdot \widehat{x} \\
&\qquad = \int_{B_{R_1}\setminus B_R} -\big(2\alpha - (d-1)\big)\big(|v|^2 - k^{-2} |\nabla v|^2\big) + k^{-2}\big(\abs{\nabla v}^2 - \abs{v_r}^2\big) + \big|k^{-1} v_r - \ri v\big|^2.\label{eq:VPS5}
\end{align}
		Setting $\alpha = (d-1)/2$ eliminates the first term on the right-hand side of \eqref{eq:VPS5}. Since $|v_r|\leq |\nabla v|$, the remaining terms on the right-hand side of \eqref{eq:VPS5} are non-negative, and thus
\beqs
\int_{\partial B_{R_1}}Q_{R_1,(d-1)/2}(u) \cdot \widehat{x} - \int_{\partial B_R} Q_{R,(d-1)/2}(u) \cdot \widehat{x} \geq 0.
\eeqs
Sending $R_1 \to \infty$ and using \eqref{eq:N0}, we obtain the result \eqref{eq:2.1}.
\epf

\bpf[Proof of Theorem \ref{thm:Morawetz}]
The plan is to integrate the identity \eqref{eq:morid1} over $B_R$ with $v=u$, $\beta=R$, and $\alpha=(d-1)/2$, and then use the divergence theorem.
We justify using the divergence theorem just as we did at the beginning of the proof of Lemma \ref{lem:MLint} to find that
if $v\in C^\infty(\overline{D})$ then 
\begin{align}\nonumber 
\int_{B_R} 2 \Re \big(\overline{\cM_{R, (d-1)/2} v } \,\cL v \big) 
= & \int_{\partial B_R} \Big[ 2 k^{-1}\Re\big(\overline{\cM_{R, (d-1)/2} v}\, k^{-1}\nabla v \big) + \big(\nvs - k^{-2}\ngvs\big) \bx\Big] \\
&-\int_{B_R}\big( \nvs +k^{-2}|\gv|^2\big).\label{eq:morid1int}
\end{align}
We now claim that  \eqref{eq:morid1int} holds for $v\in H^2(B_R)$; this follows since $C^\infty(\overline{D})$ is dense in $H^2(D)$ \cite[Page 77]{Mc:00}, and, by the trace theorem (see, e.g., \cite[Theorem 3.37]{Mc:00}), \eqref{eq:morid1int} is continuous in $v$ with respect to the topology of $H^2(B_R)$ . By Corollary \ref{cor:H2}, $u \in H^2(B_R)$, and thus \eqref{eq:morid1int} holds with $v=u$.
Using in \eqref{eq:morid1int} the definition of $Q$ \eqref{eq:Q}, the fact that $u$ is an outgoing solution of the Helmholtz equation \eqref{eq:Helmholtz}, and Lemma \ref{lem:M4}, we find that
\beqs
-\int_{B_{R}} 2\Re{\left(\overline{\mathcal{M}_{R,(d-1)/2}u}f\right)} + \norm{u}^2_{H^1_k(B_R)} = \int_{\partial B_{R}} Q_{R, (d-1)/2}(u) \cdot \widehat{x}\leq 0.
\eeqs
Thus
\beq\label{eq:VPS4}
\norm{u}^2_{H^1_k(B_R)} 
 \leq 2\norm{\mathcal{M}_{R, (d-1)/2}u}_{L^2(B_R)} \norm{f}_{L^2(B_R)}.
\eeq
By the inequality $|a + b|^2 \leq 2 |a|^2 + 2 |b|^2$, the fact that $|\ri kr + \alpha|^2 = k^2r^2 + \alpha^2$, and the bound $r \leq R$ on $B_R$, 
\beqs
\norm{\mathcal{M}_{R, (d-1)/2}u}_{L^2(B_R)}^2 \leq 2R^2 \norm{\nabla u}_{L^2(B_R)}^2 + 2[(kR)^2 + \alpha^2] \norm{u}_{L^2(B_R)}^2 \leq 2k^2R^2\left(1 + \frac{\alpha^2}{k^2 R^2}\right) \norm{u}_{H^1_k(B_R)}^2.
\eeqs
Using this in \eqref{eq:VPS4}, and recalling that $\alpha =(d-1)/2$, we find that
\beqs
\norm{u}_{H^1_k(B_R)} \leq 2kR \sqrt{1 + \left(\frac{d-1}{2kR}\right)^2} \norm{f}_{L^2(B_R)},
\eeqs
which is \eqref{eq:Morawetz_bound}.
\epf

%\bre\mythmname{Bibliographic remarks}
%The idea of multiplying second-order PDEs with first-order expressions has been used by many authors; multiplying $\Delta v$ by a derivative of $v$ goes back to Rellich \cite{Re:40, Re:43}, and multiplying $\nabla\cdot(A\gv)$ by a derivative of $v$ goes back to H\"ormander \cite{Ho:53} and Payne and Weinberger \cite{PaWe:58} (e.g., the identity \eqref{eq:morid1} with $\alpha$ and $\beta$ equal zero appears as \cite[Equation 2.4]{PaWe:58}). These identities have been independently discovered by, e.g., Jerison and Kenig \cite{JeKe:80, JeKe:81D, JeKe:81N} and Pohozaev \cite{Po:65}.
%%The Rellich identity with multiplier $\bx\cdot\gu$ is often known as the Derrick--Pohozaev identity (see, for example, \cite[\S9.4.2]{Ev:98}). This is because it was independently introduced by both Derrick \cite{De:64} and Pohozaev \cite{Po:65}) in their studies of the possible non-existence of solutions of $\Delta u + \lambda f(u)=0$
%%(although DerrickÕs proof of the identity was not via the multiplier $\bx\cdot\gu$, but via a scaling argument: see \cite[\S2.1]{BeLi:83} for a nice account of both proofs.)
%In the context of the Helmholtz equation, the identity \eqref{eq:morid1} with $\bx$ replaced by a general vector field, and $\alpha$ and $\beta$ replaced by general scalar fields was the heart of Morawetz's paper \cite{Mo:75}, following both the earlier work by Morawetz and Ludwig \cite{MoLu:68} using \eqref{eq:ml2d} and Morawetz's earlier work on the wave equation \cite{Mo:61}.
%\ere

\section*{Acknowledgements}
The idea for this paper came out of the ``Nachdiplom'' lecture course I taught at ETH Z\"urich in Fall 2021; it is a pleasure to thank the Institute for Mathematical Research (FIM) and the Seminar for Applied Mathematics at ETH for their hospitality during that time.
In particular, the proof of Theorem \ref{thm:splitting} presented in this paper was the result of Ralf Hiptmair asking me about the minimum technicalities needed for the proof  in  \cite{LSW3} of the analogous splitting in the variable-coefficient case.
I also thank Martin Averseng for many useful discussions in his role as the teaching assistant for this course. 
Finally I thank Alastair Spence (University of Bath), Jared Wunsch (Northwestern University), and an anonymous referee for their insightful comments on earlier versions of the paper, and I acknowledge support from EPSRC grant  EP/R005591/1.

%%%%%%%%%%%%%%%%%%%%%%%%%%%%%%%%%%%%%%%%%%%%%%%%%%%%%%%%%%%%%%%%%%%%%%%%%%%%%%%%%%%%%%%%%%%%%%%%%%%%%%%%%%%%

\footnotesize{
\bibliographystyle{plain}
\bibliography{biblio_acta}

\begin{thebibliography}{10}
\bibitem{Wh:15}
E.~T. Whittaker.
\newblock On the functions which are represented by the expansions of the
  interpolation-theory.
\newblock {\em Proceedings of the Royal Society of Edinburgh}, 35:181--194,
  1915.
\end{thebibliography}
}

\end{document}